\documentclass[reqno,11pt]{article}
\usepackage{amsmath,amsfonts,amssymb,euscript,theorem}

\numberwithin{equation}{section}
\newtheorem{theorem}{Theorem}[section]

{\theorembodyfont{\rmfamily} \newtheorem{remark}[theorem]{Remark} }
\newtheorem{lemma}[theorem]{Lemma}

\newcommand{\bke}[1]{\left( #1 \right)}

\newcommand{\bket}[1]{\left\{ #1 \right\}}
\newcommand{\norm}[1]{\left\Vert #1 \right\Vert}
\newcommand{\abs}[1]{\left| #1 \right|}
\newcommand{\lan}{\langle}
\newcommand{\ran}{\rangle}

\newcommand{\al}{\alpha}
\newcommand{\e}{\varepsilon}
\newcommand{\de}{\delta}

\newcommand{\ph}{\varphi}
\newcommand{\si}{\sigma}
\renewcommand{\th}{\theta}
\newcommand{\R}{\mathbb{R}}
\newcommand{\C}{\mathbb{C}}
\newcommand{\Z}{\mathbb{Z}}
\renewcommand{\Re}{\mathop{\mathrm{Re}}}

\newcommand{\pd}{\partial}
\newcommand{\donothing}[1]{{}}

\newcommand{\sech}{\,\mathrm{sech}\,}

\newcommand{\sE}{\EuScript{E}}
\newcommand{\sO}{\EuScript{O}}

\newcommand{\dist}{\, \mathrm{dist}\,} 

\newcommand{\Emin}{\sE_\text{min}}

\newcommand{\tf}{\tilde{f}}

\newenvironment{proof}{{\bf Proof.}}{\hfill\fbox{}\par\vspace{.2cm}}

\begin{document}

\title{Schr\"odinger Flow Near Harmonic Maps}

\author{{}\\[0mm]
Stephen Gustafson \quad Kyungkeun Kang \quad Tai-Peng Tsai \\
        \;\; \small{gustaf@math.ubc.ca} \quad \;\; 
\small{kkang@math.ubc.ca} \quad \;\small{ttsai@math.ubc.ca} \\[5mm]
Department of Mathematics, 
        University of British Columbia \\
Vancouver, Canada, V6T1Z2}
\date{}

\maketitle


\renewcommand{\S}{\mathbb{S}}

\begin{abstract}
  For the Schr\"odinger flow from
  $\R^2\times\R^+$ to the $2$-sphere $\S^2$, it is not known if
  finite energy solutions can blow up in finite time.
  We study equivariant solutions whose energy is near the energy 
  of the family of equivariant harmonic maps. We prove that such solutions
  remain close to the harmonic maps until the blow up time (if any), 
  and that they blow up if and only if the length scale of the nearest 
  harmonic map goes to zero.
\end{abstract}


\tableofcontents

\section{Introduction and main result}
\label{sec:intro}

The Schr\"odinger flow for maps from $\R^n$ to $\S^2$ 
(also known as the Schr\"odinger map, and, in ferromagnetism, as the
Heisenberg model or Landau-Lifshitz equation)
is given by the equation
\begin{equation}
  u_t = u \times \Delta u, \quad\quad  u(x,0)=u_0(x).
\label{int:100}
\end{equation}
Here $u(x,t)$ is the unknown map from $\R^n \times \R^+$ to the
$2$-sphere 
\[
  \S^2 := \{ u \in \R^3 \;\; | \;\; |u|=1 \} \subset \R^3,
\]
$\Delta$ denotes the Laplace operator in $\R^n$, and 
$\times$ denotes the cross product in $\R^3$.
A more geometric way to write this equation is
\begin{equation}
  u_t=J P \Delta u, \quad\quad P \Delta u = \Delta u + |\nabla u|^2 u
\end{equation}
where $P=P_u$ denotes the orthogonal projection from $\R^3$ onto the
tangent plane 
\[
  T_u \S^2 := \{ \xi \in \R^3 \; | \; \xi \cdot u = 0 \}
\]
to $\S^2$ at $u$, and $J = J^u := u \times$ is a rotation 
through $\pi/2$ on $T_u \S^2$.  

On one hand, Equation~\eqref{int:100} is a borderline case of the 
Landau-Lifshitz-Gilbert equations which model isotropic ferromagnetic 
spin systems:
\begin{equation}
  u_t = a P \Delta u + b J P \Delta u, \quad a \ge 0.
\end{equation}
(see, eg., \cite{LL,KIK}).
The Schr\"odinger flow corresponds to the case $a=0$. The case $b=0$
is the well-studied harmonic map heat flow, for which some
finite-energy solutions do blow up in finite time (\cite{cdy}).  

On the other hand, equation~\eqref{int:100} is a particular case of
the Schr\"odinger flow from a Riemannian manifold into
another one with a complex structure 
(see, eg., \cite{ZLT,DW98,TU,DW01,GS,PWW,D02}).  
We will limit ourselves to the case 
$u:\R^n \times \R^+ \to \S^2$ in this paper.

Equation \eqref{int:100} can be written in the divergence form $u_t =
\sum_{j=1}^n \pd_j (u \times \pd_j u)$, which is useful in the
construction of global weak solutions \cite{SSB}. Its formal 
equivalence to a nonlinear Schr\"odinger equation (NLS) can be seen by 
applying the stereographic projection from $\S^2$ to $\C_\infty$, the extended
complex plane:
\begin{equation}
  \label{eq:1-5}
  w = \frac {u_1+ i u_2}  {1+u_3}, \quad
  i w_t = -\Delta w + 
  \frac {2 \bar w}{1+|w|^2}{\textstyle \sum_j} (\pd_j w)^2.
\end{equation}
It is also known to be equivalent to an integrable 
cubic NLS in space dimension $n=1$
(see, eg., \cite{DPL,T}).

Equation \eqref{int:100} formally conserves the {\it energy}
\begin{equation}
  \label{eq:1-6}
  \sE(u) = \frac 12 \int_{\R^n} |\nabla u|^2 \, dx 
  = \frac 12 \int_{\R^n} \sum_{j=1}^n \sum _{k=1}^3 |\pd_j u_k|^2 \, dx.
\end{equation}
The space dimension $n=2$ is critical in the
sense that $\sE(u)$ is invariant under scaling.  In general,
\begin{equation}
  \sE(u) = s^{n-2} \sE(u_s), \quad u_{s}(x) := u(x/s), \quad s >0.
\end{equation}

For our problem $u:\R^n \times \R^+ \to \S^2$, local in time
well-posedness (LWP) is established in \cite{SSB} in the class $|u|=1$
and $\nabla u \in H^k(\R^n)$, where $k>n/2+1$ is an integer.
They also proved global in time well-posedness (GWP) in the same class
when $n=1$, and when $n\ge 2$ for data which is small in certain 
Sobolev norms.
For $n=2$, global existence is proved in \cite{CSU}
for small energy radial or equivariant data.
Also for $n=2$, LWP for a closely related system of nonlinear Schr\"odinger 
equations is established in \cite{NSU03,NSU04},
for data corresponding to $\nabla u \in H^{1+\epsilon}$.

There are known to be self-similar blow-up solutions 
for $n=2$ \cite{CSZ}; however, these do {\it not} have finite energy.

Solutions resembling solitary waves exist for $n=2$ when the target
manifold is the hyperbolic 2-space, but these, again, have infinite
energy, see \cite{D04}.  They generate blow-up solutions as well.

Fix $m \in \Z$ a non-zero integer. 
By an $m$-equivariant map $u:\R^2 \to \S^2$,
we mean a map of the form
\begin{equation}
 u(r,\th) =  e^{m \th R} \, v(r)
\end{equation}
where $(r,\th)$ are polar coordinates on $\R^2$,
$v : [0, \infty) \to \S^2$, and $R$ is the 
matrix generating rotations around the $u_3$-axis:
\begin{equation}
R=\begin{bmatrix}0&-1&0\\1&0&0\\0&0&0\end{bmatrix},
\quad
e^{\al R}=\begin{bmatrix}\cos\al&-\sin\al&0\\ 
\sin\al&\cos\al&0\\0&0&1\end{bmatrix}.
\end{equation}
Radial maps arise as the case $m=0$.
The class of $m$-equivariant maps is formally
preserved by the Schr\"odinger flow.

If $u$ is $m$-equivariant, we have 
$|\nabla u|^2 = |u_r|^2 + r^{-2} |u_\th|^2
= |v_r|^2 + \frac{m^2}{r^2}|Rv|^2$ and so
\begin{equation}
  \label{eq:1-10}
  \sE(u) = \pi\int_0^\infty 
  \bke{|v_r|^2 + \frac {m^2}{r^2} (v_1^2+v_2^2)} \, r dr .
\end{equation}
If $\sE(u)<\infty$, the limits $\lim_{r \to 0} v(r)$ and
$\lim _{r \to \infty} v(r)$ make sense (see~(\ref{eq:limits})
in the next section), and so we must have 
$v(0),v(\infty) = \pm \hat{k}$, where $\hat{k} = (0,0,1)^T$.
We may and will fix $v(0)=-\hat{k}$. The two cases 
$v(\infty) = \pm \hat{k}$ correspond to different topological 
classes of maps.  
We denote by $\Sigma_m$ the class of $m$-equivariant maps with 
$v(\infty) = \hat{k}$:
\begin{equation}
  \Sigma_m = \bket{u: \R^2 \to \S^2 \; | \quad u=e^{m \th R}v(r), \, 
  \sE(u)<\infty, \, v(0)=-\hat{k},\, v(\infty)=\hat{k}}.
\end{equation}

The energy $\sE(u)$ can be rewritten as follows:
\begin{equation}
\label{eq:1-11}
  \sE(u) = \pi \int_0^\infty \bke{|v_r|^2 + 
  \frac {m^2}{r^2} |J^v  R v|^2} \, r dr
  = \pi \int_0^\infty |v_r - \frac{|m|}{r} J^v R v|^2 \, r dr  
  + \Emin
\end{equation}
where $J^v := v \times$, and
\begin{equation}
  \Emin = 2\pi \int_0^\infty v_r \cdot \frac{|m|}{r} J^v R v \, r dr 
  = 2 \pi |m| \int_0^\infty (v_3)_r dr = 2 \pi |m| [v_3(\infty) + 1]
\end{equation}
(using $v_1^2+v_2^2+v_3^2=1$).
The number $\Emin$, which depends only on the boundary 
conditions, is in fact $4\pi$ times the absolute value of the
{\it degree} of the map $u$, considered as a map from $\S^2$
to itself (defined, for example, by integrating the pullback
by $u$ of the volume form on $\S^2$).
It provides a lower bound for the energy
of an $m$-equivariant map, $\sE(u) \geq \Emin$,
and this lower bound is attained if and only if   
\begin{equation}
  \label{eq:1-13}  v_r = \frac{|m|}{r} J^v R v.
\end{equation}
If $v(\infty) = -\hat{k}$, the minimal energy is
$\Emin=0$ and is attained by the constant map, 
$u \equiv -\hat{k}$. On the other hand, if $v(\infty) = \hat{k}$,  
the minimal energy is 
\[
  \sE(u) \geq \Emin = 4 \pi |m|
\] 
and is attained by the $2$-parameter family of harmonic maps
\begin{equation}
\label{eq:hm}
  \sO_m := \bket{e^{m \th R}h^{s,\al}(r) \; 
  | \; s > 0, \; \al \in [0, 2\pi)}
\end{equation}
where 
\[
  h^{s,\al}(r) := e^{\al R} h(r/s), 
\]
and
\begin{equation}
\label{eq:harmonic}
  h(r) = \left( \begin{array}{c} h_1(r) \\ 0 \\ h_3(r) 
  \end{array} \right), \quad
  h_1(r)=\frac 2{r^{|m|} + r^{-|m|}}, \quad
  h_3(r)= \frac {r^{|m|} - r^{-|m|}}{r^{|m|} + r^{-|m|}}.
\end{equation}
The fact that $h(r)$ satisfies~\eqref{eq:1-13} means
\begin{equation}
  (h_1)_r = -\frac m r h_1 h_3, \quad (h_3)_r = \frac m r h_1^2.
\end{equation}
So $\sO_m$ is the orbit of the single harmonic map
$e^{m \theta R} h(r)$ under the symmetries of the energy
$\sE$ which preserve equivariance: scaling, and rotation. 
Explicitly,
\begin{equation}
  e^{m \th R}h^{s,\al}(r) = \left(
  \begin{array}{c} \cos (m \th + \al)h_1( r/s) \\
  \sin(m \th + \al) h_1(r/s) \\ 
   h_3( r/s) \end{array} \right).  
\end{equation}
The solution \eqref{eq:harmonic} is easily found 
by solving the system \eqref{eq:1-13} of ODEs directly.
Alternately, under the stereographic projection \eqref{eq:1-5},
Equation~(\ref{eq:1-13}) amounts to the Cauchy-Riemann equations,
and these harmonic maps correspond to the anti-meromorphic 
(if $m > 0$) functions
\begin{equation}
  w =  e^{i \al} \left( \frac{\bar z}{s} \right)^{-m}, 
  \quad z = r e^{i \th}.
\end{equation}

We are now ready to state our main result. 
We denote $\norm{u}_{\dot H^k} = \norm{\nabla^k u}_{L^2}$.
\begin{theorem}
\label{mainthm:20}
There exist $\delta > 0$
and $C_0, C_1>0$ such that if 
$u \in C([0,T) ; \dot{H}^2 \cap \Sigma_m)$ is
a solution of the Schr\"odinger flow \eqref{int:100} 
conserving energy, and satisfying
\[
  \delta_1^2 := \sE(u_0) - 4\pi |m| < \delta^2,
\]
then there exist $s(t)\in{\mathcal C}([0,T);(0,\infty))$ and 
$\alpha(t)\in{\mathcal C}([0,T);\R)$ so that
\begin{equation}
\label{eq:mainthm}
  \norm{u(x,t)-e^{(m\theta+\alpha(t))R}h(r/s(t))}_{\dot{H}_1(\R^2)} 
  < C_0 \delta_1, \quad \forall t \in [0,T).
\end{equation}
Moreover, $s(t) > C_1/\|u(t)\|_{\dot{H}_2}$.
Furthermore, if $T < \infty$ is the maximal time of existence 
for $u$ in $\dot{H}^2$ 
(i.e. $\lim_{t \to T^-} \| u(t) \|_{\dot H^2} = \infty$), then 
\begin{equation}
\label{eq:concentrate}
  \liminf_{t\to T^{-}} s(t) = 0.
\end{equation}
\end{theorem}
\begin{remark}
\begin{enumerate}
\item
This theorem can be viewed, on one hand, as an 
{\it orbital stability} result for the family
of harmonic maps (at least up to the possible blow-up
time), and on the other hand as a characterization
of blow-up for energy near $\Emin$:
solutions blow-up if and only if the
$\dot{H}^1$-nearest harmonic map ``collapses''
(i.e. its length-scale goes to zero).  
\item
The assumption $\sE(u_0)-4\pi |m| <\delta^2$ implies $u_0$ is close to
$\sO_m$ in $\dot H^1$ (see~(\ref{eq:global})), 
but not necessarily in $\dot H^2$.
\item
The existence of local (in time) $\dot H^2 \cap \dot H^1$ solutions
of~(\ref{int:100}) is established in~\cite{SSB} for
sufficiently regular initial data. In particular, this
ensures Theorem~\ref{mainthm:20} is non-empty
(see also~\cite{NSU03}-\cite{NSU04} for local well-posedness results).
Local well-posedness with data in $\dot H^2 \cap \dot H^1$
appears still to be open. If we had this, (\ref{eq:concentrate})
could be replaced by $\lim_{t \to T^-} s(t) = 0$. 

\item From now on we will assume $m>0$.  The cases $m<0$ of Theorem
\ref{mainthm:20} follow from the change of variables $(x_1,x_2,x_3)
\to (x_1,-x_2,x_3)$.
\end{enumerate}
\end{remark}

The plan for the paper is as follows:
In Section~\ref{sec:maps} we study maps 
whose energy is close to that of the family of harmonic maps.
The analysis here is completely time-independent.
In Section~\ref{sec:flow} we apply Strichartz estimates to a certain 
nonlinear Schr\"odinger equation, obtained via the Hasimoto 
transformation introduced in \cite{CSU},
and present the proof of the main theorem.
The proofs of some of the more technical lemmas 
are relegated to Section~\ref{sec:lemmas}
in order to streamline the presentation.
Without loss of generality, we assume $m > 0$ 
for the rest of the paper.

\medskip

\noindent {\bf Remark on notation:}
throughout the paper, the letter 
$C$ is used to denote a generic constant, 
the value of which may change from line to line.


\section{Maps with energy near the harmonic map energy}
\label{sec:maps}

This section is devoted entirely to static
$m$-equivariant maps (i.e. there is no time-dependence 
anywhere in this section).
We establish some properties of maps
with energy close to the harmonic map energy $4 \pi m$.
Roughly speaking, we prove that such maps are $\dot H^1$-close
to harmonic maps. Precise statements appear in Theorem~\ref{thm:maps}
below.

We define the distance from any map $u$ to the family 
$\sO_m$ of $m$-equivariant harmonic maps to be
\[
  \dist(u,\sO_m) := \inf_{s \in (0,\infty), \al \in \S^1}
  \| u - e^{m \theta R}h^{s,\al} \|_{\dot{H}^1}.
\] 
Here $S^1 = \R/2\pi$.

The following theorem defines a (nonlinear) projection from 
the set $\Sigma_m$
of $m$-equivariant maps with energy close to $4 \pi m$
onto the family $\sO_m$, and establishes a key fact:
for maps in this set, the squared distance 
$\dist^2(u,\sO_m)$ is bounded by the energy
difference $\sE(u) - 4\pi m$. 

\begin{theorem} 
\label{thm:maps}
There are constants $\de>0$ and $C_0,C_1>0$ such that
if $u \in \Sigma_m$ satisfies 
\[
  \sE (u) < 4 \pi m + \de^2,
\] 
then the following hold:
\begin{itemize}
\item[(a)]
There exist unique $s(u) \in (0,\infty)$ and 
$\al(u) \in \S^1$ such that
\begin{equation*}
 \dist(u,\sO_m) = 
 \| u - e^{m\theta R}h^{s(u),\al(u)} \|_{\dot{H}^1}.
\end{equation*}
Moreover, $s(u)$ and $\al(u)$ are continuous functions
of $u \in \dot{H}^1$.

\item[(b)] $\dist(u,\sO_m) < C_0[\sE(u)-4\pi m]^{1/2} < C_0 \de$. 

\item[(c)] If $u \in \dot H^2(\R^2)$, 
then $s(u) \norm{u}_{\dot H^2(\R^2)} > C_1$.

\end{itemize}
\end{theorem}

\begin{proof}\,\,\,
The proof is long, so we break it into a series of steps.
At each step, we may need to take $\delta$ smaller than
in the previous step.

\medskip

{\bf Step 1: a change of variable}.
Recall that $u \in \Sigma_m$ implies that in polar coordinates $(r,\theta)$,  
$u(x) = e^{m \theta R} v(r)$, with $v(0) = -\hat{k}$ and $v(\infty) = \hat{k}$.
The change of variables
\[
  r \to y = m\log(r) \in (-\infty, \infty), 
  \quad \mbox{ or } \quad   
  e^y = r^{m}
\]
turns out to be very useful for our purposes.  Set 
\[
  \tilde{v}(y) := v(e^{y/m}). 
\]
Under this change of variables, 
the $\dot H^1$ inner-product of $m$-equivariant maps changes as follows:
\[
\begin{split}
  \lan e^{m \theta R} v(r), e^{m \theta R} w(r) \ran_{\dot H^1} 
  &= \int_{\R^2} \nabla [e^{m \theta R} v(r)] \cdot
  \nabla [e^{m \theta R} w(r)] dx \\
  &= 2\pi \int_0^\infty (v_r \cdot w_r + \frac{m^2}{r^2} Rv \cdot Rw) rdr \\
 &= 2\pi m\int_\R 
\bke{\tilde{v}' \cdot \tilde w' +  R \tilde v \cdot R \tilde w} dy,
\end{split}
\]
where ``$ \; ' \; $'' denotes $d/dy$. In particular, 
\begin{equation}\label{Etilde}
  \sE(u) = 2\pi m \tilde{E}(\tilde{v}), \quad 
  \tilde{E}(\tilde{v}) := 
  \frac{1}{2}\int_{\R}\big(\abs{\tilde{v}'(y)}^2
  +\tilde{v}_1(y)^2 + \tilde{v}_2(y)^2 \big)dy.
\end{equation}
Note that this implies $\tilde{v}_j \in H^1(\R)$
for $j = 1,2$, and in particular $(v_j^2)' \in L^1(\R)$
so that the limits $lim_{y \to \pm \infty} v_j^2(y)$
exist, and are equal $0$. $\tilde{v}$ is continuous,
and $\tilde{v}_3^2 = 1 - \tilde{v}_1^2 - \tilde{v}_2^2$
has limit $1$ as $y \to \pm \infty$. Thus
the limits $\lim_{y \to \pm \infty} \tilde{v}$ exist,
justifying our earlier claim 
\begin{equation}
\label{eq:limits}
  \lim_{r \to 0} v(r) \quad \mbox{ and } \quad
  \lim_{r \to \infty} v(r) \quad \mbox{ exist. }
\end{equation}
Recall that for $v \in \Sigma_m$, we have chosen
$\tilde{v}(-\infty) = -\hat{k}$, $\tilde{v}(\infty) = \hat{k}$.

$\tilde{E}(\tilde{v})$ inherits the ``topological lower bound''
\begin{equation}
\label{eq:topE}
  \tilde{E}(\tilde{v}) = 2 + 
  \frac{1}{2} \int_{\R}|\tilde{v}' - J^{\tilde{v}} R \tilde{v}|^2 dy   
  \geq 2.
\end{equation}

In this new variable, scaling $r$ corresponds to translating $y$.
In particular, the family of harmonic maps is composed 
of translations and rotations of a simple explicit map: 
\[
\begin{split}
  &h^{s,\al}(r) = e^{\al R} h(r/s) =   
  e^{\al R} \tilde{h}(y-m\log(s)) \\
  & \tilde{h}(y) = \left( 
  \begin{array}{c} \sech y \\ 0 \\ \tanh y \\
  \end{array} \right).
\end{split}
\]

\medskip

{\bf Step 2: energy close to $4\pi m$ implies $u$ close to a harmonic map}.
Expressed in the variable $y$, what we would like to prove is
\begin{lemma}
\label{Lemma3-3}
For any $\epsilon > 0$, there exists $\mu > 0$ such that if
a map $\tilde{v} : \R \to \S^2$ satisfies 
$\tilde{E}(\tilde{v}) < 2 + \mu$, and 
$\tilde{v}(-\infty) = -\hat{k}$, $\tilde{v}(\infty) = \hat{k}$,
then
\[
  \inf_{\al\in{\mathbb{S}^1},a\in\R}
  \norm{\tilde{v}-e^{\al R}\tilde{h}(\cdot-a)}_{H^1(\R)}<\epsilon.
\] 
\end{lemma}

\begin{proof}
Suppose not. Then there exist $v_j(r)$, $j = 1,2,3,\ldots$, 
and $\epsilon_0>0$ such that 
\begin{equation}
  |v_j(r)| \equiv 1, \;
  \tilde{E}(v_j)<2+\frac{1}{j},\,\, 
  v_{j}(\pm \infty) = \pm \hat{k},\,\,
  \norm{v_j-e^{\al R}\tilde{h}(\cdot-a)}_{H^1(\R)}\geq\epsilon_0
\label{lemma3-3:10}
\end{equation}
for every $\al\in\mathbb{S}^1$ and $a\in\R$.
Since $\tilde{E}(v_j)<\infty$, $v_j$ is continuous, and thus 
$v_{j3}(a_j)=0$ for some $a_j\in\R$.
We replace $v_j$ by $w_j(y):=e^{-\al_jR}v_j(y+a_j)$, where
$\al_j\in{\mathbb{S}^1}$ is chosen so that $w_{j}(0)=\hat{\imath}=(1,0,0)^T$.
The properties \eqref{lemma3-3:10} still hold for $\{w_j\}$:
\begin{equation}
  |w_j(r)| \equiv 1,\; 
  \tilde{E}(w_j)<2+\frac{1}{j}, \,\, 
  w_{j}(\pm \infty) = \pm \hat{k},\,\,
  \norm{w_j-\tilde{h}}_{H^1(\R)}\geq\epsilon_0.
\label{lemma3-3:20}
\end{equation}
Since $\sup_j \tilde{E}(w_j)<\infty$, there is a 
subsequence (which we continue to denote by $\{w_j\}$),
and a limit vector function $w^*(y)$, satisfying
\begin{eqnarray}
\left\{
\begin{array}{cc}
  w_{j1} \,\, \longrightarrow \,\, w^*_{1}, \quad 
  w_{j2} \,\, \longrightarrow \,\, w^*_{2}
  &\mbox{ weakly in }H^1(\R),\\
  w'_{j3} \,\,\longrightarrow \,\, w^{*'}_3 &\mbox{ weakly in }L^2(\R),\\
  w_{j}\,\,\longrightarrow\,\, w^*  
  &\mbox{ strongly in }L^2_{\rm{loc}}(\R)\cap
  {\mathcal{C}}^{0}_{\rm{loc}}(\R). 
\end{array} 
\right.
\end{eqnarray}
Because of the local uniform convergence,
we have $|w^*| \equiv 1$ and $w^*(0)=\hat{\imath}$. 
On the other hand, by the topological lower 
bound~(\ref{eq:topE}), we have
\[
  2 \leq \tilde{E}(w_j) = 2 + \frac{1}{2}\int_{\R}
  \abs{w'_j-J^{w_j}Rw_j}^2<2+\frac{1}{j},
\]  
from which it is immediate that $\norm{w'_j-J^{w_j}Rw_j}_{L^2(\R)}
\rightarrow 0$ as $j\rightarrow\infty$.
For any bounded interval $I$, using 
$J^w Rw = \hat{k} - w_3w$, we have
\[
  \norm{J^{w_j}Rw_j-J^{w^*}Rw^*}_{L^2(I)} =
  \norm{w_{j3}w_j - w^*_3 w^*}_{L^2(I)}
  \to 0 \mbox{ as } j \to \infty.
\]
Hence $w_j' \to J^{w^*} R w^*$ strongly in $L^2(I)$.
On the other hand, ${w^*}'$ is the weak limit of $w'_j$, 
and so we obtain ${w^*}'=J^{w^*}Rw^*$ almost everywhere,
with $w^*(0)=\hat{\imath}$. By uniqueness of $H^1_{loc}$ solutions
of this system of ordinary differential equations, 
we must have
\[
  w^*(y) = \tilde{h}(y) = (\sech y, 0 , \tanh y)^T.
\]
Now we note that 
\[
  \tilde{E}(w^*) = \tilde{E}(\tilde{h}) = 2 =
  \lim_{j\rightarrow \infty} \tilde{E}(w_j).
\]
By \eqref{Etilde},
\[
  \norm{Rw^*}^2_{H^1(\R)} + \norm{{w^*_3}'}^2_{L^2(\R)}=
  \lim_{j\rightarrow \infty}\big(\norm{Rw_j}^2_{H^1(\R)}
  +\norm{w'_{j3}}^2_{L^2(\R)}\big).
\]
By weak lower semi-continuity, i.e.
\[
  \norm{Rw^*}^2_{H^1(\R)}\leq \liminf_{j\rightarrow 0}\norm{Rw_j}^2_{H^1(\R)},
  \quad \norm{w^{*'}_3}^2_{L^2(\R)}\leq \liminf_{j\rightarrow 0}
  \norm{w'_{j3}}^2_{L^2(\R)},
\]
we have 
$\norm{Rw_j}^2_{H^1(\R)}\rightarrow \norm{Rw^*}^2_{H^1(\R)}$ and 
$\norm{w'_{j3}}^2_{L^2(\R)}\rightarrow  \norm{w^{*'}_3}^2_{L^2(\R)}$,
which implies $Rw_j\rightarrow Rw^*$ strongly in $H^1(\R)$ and
$w'_{j3}\rightarrow w^{*'}_3$ strongly in $L^2(\R)$.
Finally, we will show that $w_j-\tilde{h}$ converges to $0$ strongly 
in $H^1(\R)$, which will contradict assumption \eqref{lemma3-3:20},
and so complete the proof of the lemma. 
Indeed,
\[
  \norm{w_j-\tilde{h}}^2_{H^1(\R)} = \norm{w'_j-\tilde{h}'}^2_{L^2(\R)}
  +\norm{Rw_j-R\tilde{h}}^2_{L^2(\R)}+\norm{w_{j3}-\tilde{h}_{3}}_{L^2(\R)}^2.
\]
We have already shown that the first two terms go to zero in the limit,
and so it remains to consider the last term.
For this, we need another lemma.
For $f:(a,b)\rightarrow \R$, denote by $T_{(a,b)}(f)$ 
the total variation of $f$ on $(a,b)$.
The following lemma shows that the total variation of 
$v_3$ is close to $2$ if $\sE(u)$ is close to $4\pi m$.
\begin{lemma}
\label{Lemma3-4}
  If $u = e^{m \theta R} v(r) \in \Sigma_m$ and 
  $\sE(u) = 4 \pi m + \epsilon_0$, then
  $T_{(0,\infty)}(v_3) \leq 2+C\epsilon_0$.
\end{lemma}
\begin{proof}
Make the change of variable
$\tilde{v}(y) = v(e^{y/m})$, and write
$\tilde{v} = (\tilde{\rho} \cos(\tilde{\omega}), 
\tilde{\rho} \sin(\tilde{\omega}), \tilde{v}_3)$,
so that $\tilde{\rho}^2 = \tilde{v}_1^2 + \tilde{v}_2^2$.
We have
\[
\begin{split}
  4+\frac{\epsilon_0}{\pi m} &\geq 
  \int_{\R}\big(\abs{\tilde{v}'}^2
  +\abs{R \tilde{v}}^2\big) dy 
  =\int_{\R}\big(\tilde{\rho}^2 (\tilde{\omega}')^2 
  +\frac{(\tilde{\rho}')^2}{1-\tilde{\rho}^2} + \tilde{\rho}^2\big)dy \\
  &\geq 2\int_{\R}\frac{\tilde{\rho}\abs{\tilde{\rho}'}}
  {(1-\tilde{\rho}^2)^{1/2}}dy
  =2\int_{\R}\abs{\frac{d}{dy}(1-\tilde{\rho}^2)^{\frac{1}{2}}}dy
  =2T_{\R}(\tilde{v}_3).
\end{split}
\]
Dividing by $2$ on both sides completes the proof.
\end{proof}

Applying this lemma to $w_j$, we have
$T_{\R}(w_{j3}) \leq 2 + C/j$.
Since $w_{j3}(-\infty) = -1$, $w_{j3}(0) = 0$, and
$w_{j3}(\infty) = 1$, we have $w_{j3}(y) > - C/j$
for $y \geq 0$ (and similar for $y \leq 0$). 
Fix $\epsilon_1 > 0$. For $|y| \geq \epsilon_1$ and 
 $j$ sufficiently large (depending on $\epsilon_1$),
\[
\begin{split}
  |w_{3j}-\tilde{h}_3| &=
  \frac{|w_{3j}^2-\tilde{h}_3^2|}{|w_{3j}+\tilde{h}_3|}
  \leq \frac{|w_{3j}^2-\tilde{h}_3^2|}{|\tilde{h}_3| - C/j} \\
  &\leq \frac 2{\tanh \epsilon_1}  |w_{3j}^2-\tilde{h}_3^2|
  = \frac 2{\tanh \epsilon_1} ||R\tilde{h}|^2 - |Rw_j|^2|,
\end{split}
\]
and so
\[
  \int_{|y| \geq \epsilon_1} |w_{3j}-\tilde{h}_3|^2 dy
  \to 0 \mbox{ as } j \to \infty.
\]
Since 
$\int_{|y| < \epsilon_1} |w_{3j}-\tilde{h}_3|^2 dy \leq C \epsilon_1$
and $\epsilon_1$ is arbitrary, we conclude
$\| w_{3j} - \tilde{h}_3 \|_{L^2(\R)} \to 0$.
This completes the proof of Lemma~\ref{Lemma3-3}.
\end{proof}

Translating Lemma~\ref{Lemma3-3} back to the original
variable $r = e^{y/m}$, we find:
\begin{equation}
\label{eq:global}
\begin{split}
  &\mbox{given } \epsilon > 0, \mbox{ there is }
  \mu > 0 \mbox{ s.t. if } \\
  &\mbox{      } u \in \Sigma_m, \; \sE(u) < 4\pi m + \mu, \\
  &\mbox{then } \dist(u,\sO_m) < \epsilon.
\end{split}
\end{equation}

\medskip

{\bf Step 3: existence of $s(u)$ and $\al(u)$}.
Recall that since $u \in \Sigma_m$, 
$\sE(u) \geq 4\pi m$, and 
set $\delta_1 := [\sE(u)-4\pi m]^{1/2} < \de$.
We observe first that 
\begin{equation}
  \lim_{s \to \infty} \inf_\al 
  \| u - e^{m\theta R} h^{s,\al} \|^2_{\dot H^1}
  =\lim_{s \to 0} \inf_\al \| u - e^{m\theta R} h^{s,\al} \|^2_{\dot H^1} 
  = 8\pi m + \delta^2_1.
\label{Prop3-5:10}
\end{equation}
Indeed, we have 
\begin{equation}
   \| u - e^{m\theta R} h^{s,\al} \|^2_{\dot H^1}
   =8\pi m + \delta^2_1 - 2\int_{\R^2} 
   \nabla u\cdot\nabla (e^{m\theta R}h^{s,\al})dx,
\label{Prop3-5:20}
\end{equation}
and it suffices to show that for any 
$\alpha\in\mathbb{S}^1$ 
\begin{equation}
  \int_{\R^2}\nabla u\cdot\nabla (e^{m\theta R}h^{s,\al}) dx
  \rightarrow 0,
\label{Prop3-5:30}
\end{equation}
as $s\rightarrow 0$ or $s\rightarrow\infty$.
Since $h^{s,\al}(r) = h^{1,\al}(r/s)$, this latter fact
follows from an easy lemma:
\begin{lemma}
\label{lem:scaling}
If $f \in L^2(\R^2)$, then
\[
  \frac{1}{s} f(x/s) \to 0 \mbox{ weakly in } L^2(\R^2)
  \mbox{ as } s \to \infty \mbox{ and } s \to 0.
\]
\end{lemma}
\begin{proof}
First suppose $f \in L^2 \cap L^\infty$,
and fix $g \in L^2$. By H\"older's inequality,
\[
\begin{split}
  \int_{\R^2} \frac{1}{s} f(x/s) g(x) dx
  &= \int_{|x| \leq \sqrt{s}} \frac{1}{s} f(x/s) g(x) dx
  + \int_{|x| > \sqrt{s}}  \frac{1}{s} f(x/s) g(x) dx \\
  & \leq \| f \|_{L^\infty} \| g \|_{L^2} \frac{\sqrt{\pi}}{\sqrt{s}}
  + \| f \|_{L^2} \left( \int_{|x|>\sqrt{s}} |g(x)|^2 dx \right)^{1/2}
  \to 0
\end{split}
\]
as $s \to \infty$. A similar argument covers the $s \to 0$ case.
For general $f \in L^2$,
choose $f_\epsilon \in L^2 \cap L^\infty$ with
$\| f_\epsilon - f \|_{L^2} < \epsilon$.  Then
by the above argument,
\[
  |\int_{\R^2} \frac{1}{s} f(x/s) g(x) dx|
  \leq \epsilon \| g \|_{L^2} + o(s).
\] 
Since $\epsilon$ is arbitrary, we are done.
\end{proof} 

To prove the claim~(\ref{Prop3-5:30}),
just take $g = \nabla u$ and 
$\frac{1}{s}f(\cdot/s) = \nabla (e^{m\theta R}h^{s,\al})$
in this Lemma. 

So on one hand,
$\inf_{\al} \| u - e^{m \theta R} h^{s,\al} \|^2_{\dot H^1}$
approaches $8 \pi m + \delta_1^2$ as $s \to 0$ and $s \to \infty$.
On the other hand,~(\ref{eq:global}) shows that
if $\delta$ (and hence $\delta_1$) is sufficiently small,
then for some $s$ and $\alpha$, 
$\| u - e^{m \theta R} h^{s,\al} \|^2_{\dot H^1} < 8\pi m$.
Thus to minimize
\[
  F(s,\al) := \| u - e^{m \theta R} h^{s,\al} \|^2_{\dot H^1},
\]
over $s \in (0,\infty)$, $\al \in \S^1$,
it suffices to consider $s$ in a compact subset of $(0,\infty)$.
Since $F(s,\al)$ is continuous, there must exist 
$s(u) \in (0,\infty)$ and $\al(u) \in \mathbb{S}^1$ such that 
\[
  \dist (u,\sO_m) = 
  \| u - e^{m\theta R}h^{s(u),\al(u)} \|_{\dot H^1}.
\]

\medskip

{\bf Step 4: uniqueness of $s(u)$ and $\al(u)$}.
Denote $\si = (s,\al)$. 
Suppose there exist $\si_1, \si_2$ with $\si_1\neq \si_2$ such that
\[
  \de_0: = \dist (u,\sO_m) 
  = \| u - e^{m\theta R}h^{\si_1} \|_{\dot H^1} 
  = \| u - e^{m\theta R}h^{\si_2} \|_{\dot H^1}.
\]
Let $\mu$ be half the distance between $e^{m\theta R}h^{\si_1}$ 
and $e^{m\theta R}h^{\si_2}$:
$\mu := (1/2)\| e^{m\theta R}h^{\si_1}-e^{m\theta R}h^{\si_2} \|_{\dot H^1}$. 
It follows that $\mu \le  \de_0$.
Now set $\ph(t) = e^{m\theta R}h^{\si(t)}$ with 
$\si(t) = \frac 12[(\si_1+\si_2)+t(\si_2-\si_1)]$, so that 
$\ph(-1) = e^{m\theta R}h^{\si_1}$ and 
$\ph(1) = e^{m\theta R}h^{\si_2}$.
Set $\bar{\ph} := \frac 12[ \ph(-1) + \ph(1) ]$.
Lemma~\ref{lem:curvature} (stated and proved in Section~\ref{sec:lemmas}) 
then yields
\[
  \| \bar{\ph} - \ph(0) \|_{\dot H^1}
  < C \mu^2.
\]
This estimate amounts to a bound on the curvature
of the family $\sO_m$.  

Now $\nabla[u-\bar \ph]$ is $L^2$-orthogonal to 
$\nabla[e^{m\theta R}(h^{\si_1} -h^{\si_2})]$, 
since $\|u - e^{m\theta R}h^{\si_1}\|_{\dot H^1} =
\| u - e^{m\theta R}h^{\si_2} \|_{\dot H^1}$. 
By the Pythagorean theorem,
\[
  \| u - \bar \ph \|_{\dot H^1}^2 
  = \| u - \ph(-1) \|_{\dot H^1}^2  
  - \frac 14\, \| \ph(-1) - \ph(1) \|_{\dot H^1}^2 
  = \delta^2_0 - \mu^2,
\]
and so
\[
\begin{split}
  \| u - \ph(0) \|_{\dot H^1} &\le 
  \| u - \bar\ph \|_{\dot H^1} + \| \bar\ph - \ph(0) \|_{\dot H^1}
  \le (\de^2_0- \mu^2)^{1/2} + C \mu^2 \\
  &= \de_0 - \frac {\mu^2}{\de_0 + \sqrt{\de^2_0-\mu^2}} + C \mu^2.
\end{split}
\]
By~(\ref{eq:global}), we can ensure $\de_0 < 1/(2C)$
by choosing $\delta$ sufficiently small.
This, in turn, implies
$\| u - \ph(0) \|_{\dot H^1} < \delta_0$,
which contradicts the assumption that $\de_0$ is the minimal distance.
This establishes uniqueness of $s(u)$ and $\al(u)$.

\medskip

{\bf Step 5: continuity of $s(u)$ and $\al(u)$}. 
We could invoke the implicit function theorem,
but we prefer to give a simple direct proof of continuity.
Suppose $u_j \to u$ in $\dot H^1$ with
$\sE(u_j) < 4\pi m + \delta^2$ and
$\sE(u) < 4\pi m + \delta^2$.
We have
\[
  \dist(u_j,\sO_m ) 
  \leq \| u_j - e^{m \theta R}h^{s(u),\al(u)} \|_{\dot H^1}
  \leq \| u_j - u \|_{\dot H^1} 
  + \dist(u,\sO_m)
\]
and
\[
  \dist(u,\sO_m ) 
  \leq \| u - e^{m \theta R} h^{s(u_j),\al(u_j)} \|_{\dot H^1}
  \leq \| u_j - u \|_{\dot H^1} 
  + \dist(u_j,\sO_m)
\]
and so $\dist(u_j,\sO_m) \to \dist(u_*,\sO_m)$.
Since
\[
  \| e^{m \theta R}(h^{s(u_j),\alpha(u_j)} - h^{s(u),\alpha(u)}) \|_{\dot H^1}
  \leq \dist(u_j,\sO_m) + \dist(u,\sO_m) +
  \|u_j - u\|_{\dot H^1},
\]
 $\{ s(u_j) \}$ is contained
in a compact subinterval of $(0,\infty)$ by Lemma \ref{Lemma2-3}, and so,
up to subsequence, $s(u_j) \to s_*$ and $\al(u_j) \to \al_*$,
for some $s_*$ and $\al_*$. Along this subsequence
\[
\begin{split}
  \norm{u - h^{s(u),\al(u)}}_{\dot{H}^1}
  &=\dist(u,\sO_m)
  =\lim_{j\rightarrow\infty} \dist(u_j,\sO_m) \\
  &=\lim_{j\rightarrow\infty}\norm{u_j-h^{s(u_j),\al(u_j)}}_{\dot{H}^1}
  =\norm{u - h^{s_*,\al_*}}_{\dot{H}^1}.
\end{split}
\]
By the uniqueness we have already proved, 
$s_* = s(u)$ and $\al_* = \al(u)$.
We conclude that $s(u_j) \to s(u)$, and $\al(u_j) \to \al(u)$
(for the full sequence). Continuity is proved. 
This completes the proof of part (a)
of Theorem~\ref{thm:maps}.

\medskip

{\bf Step 6: the ``linearized operator''.}
We now proceed to the proof of part (b)
of Theorem~\ref{thm:maps}. The main idea is this:
the ``global'' result of Lemma~\ref{Lemma3-3} allows us 
now to work ``locally'' -- i.e. nearby a harmonic map. 
Indeed to prove (b), we study the second 
variation of the energy functional
around the nearby harmonic map. We begin by discussing
this ``linearized operator''.

Given an $m$-equivariant map $u \in \Sigma_m$ with
$\sE(u) - 4\pi m < \delta^2$, we fix 
$s = s(u)$ and $\alpha = \alpha(u)$ and write 
$u = e^{m\th R}v(r)$ with 
\begin{equation} \label{xi-def}
  v(r) = e^{\al R}(h(r/s) + \xi(r/s)).
\end{equation}
This defines $ \xi=\xi(r) = e^{-\al R} v(sr) - h(r)$
(note that the variable $r$ here is no longer the original polar coordinate). 
Using $h_r - \frac mr J^hRh=0$, expand
\begin{align}
\nonumber
  \sE(u) &= \frac 12 \int |\nabla u|^2 dx
  =4 \pi m + \frac{1}{2}\int_{\R^2} |v_r - \frac{m}{r} J^vRv|^2 dx \\
  &= 4\pi m + \frac{1}{2}\int_{\R^2}|L \xi + \frac{m}{r} \xi_3 \xi|^2 dx
\label{Eu-expand}
\end{align}
where $L \xi$ is the linear part,
\begin{equation}
  L \xi := \xi_r + \frac mr ( \xi_3 h + h_3 \xi ).
\label{lin-oper}
\end{equation}
In fact, the operator $L$ maps tangent vector fields
(i.e. vector functions $\eta(r)$ tangent to $\S^2$ at $h(r)$), 
into tangent vector fields.
To see this explicitly, we specify an orthonormal basis
of $T_h \S^2$:
\[
  e = \left( \begin{array}{c} 0 \\ 1 \\ 0 \end{array} \right)
  \quad \mbox{ and } \quad
  J^h e = \left( \begin{array}{c} -h_3(r) \\ 0 \\ h_1(r)  \end{array} \right).
\]
Then for any map $\xi : [0,\infty) \to \R^3$, 
we have an orthogonal decomposition
\begin{equation} \label{xi-decomp}
  \xi(r) = z_1(r) e + z_2(r) J^h e + \gamma(r) h,
\end{equation}
which defines a complex-valued function 
$z(r) := z_1(r) + iz_2(r)$.  Note that
\[
L e = \frac mr h_3e, \quad 
L(J^h e) = \frac mr h_3 J^h e, \quad 
L h = \frac mr h_3 h + \frac mr \hat k.
\]
Hence the operator $L$ restricted to $T_h \S^2$ is equivalent to
\[
  L_0 := \pd_r + \frac mr h_3 = h_1 \pd_r \frac{1}{h_1},
\]
in the sense that
\[
  L (z_1 e + z_2 J^h e )
= (L_0 z_1)e + (L_0 z_2)J^h e.
\]

Roughly speaking, our strategy for proving part (b)
of Theorem~\ref{thm:maps} is to show that 
$\dist (u,\sO_m )$
is controlled by $z$, which is controlled by
$L_0 z$, which in turn is controlled by $\sE(u) - 4\pi m$. 

Define for radial complex-valued functions
$f(r)$ and $g(r)$ the following inner-product:
\[
  \lan f, g \ran_X := 
  \int_0^\infty \left(
  \bar{f}_r(r) g_r(r) + \frac{m^2}{r^2} \bar{f}(r) g(r)
  \right) rdr.
\]
If we set
$\tilde{f}(y) := f(e^{y/m})$ and 
$\tilde{g}(y) := g(e^{y/m})$, we have
\begin{equation}
\label{eq:inner}
 \lan f, g \ran_X = 
   m \int_\R (\overline{\tilde{f}'} \tilde{g}' 
+ \overline{\tilde{f}} \tilde{g}) dy =
   m \langle \tilde{f}, \tilde{g} \rangle_{H^1(\R)}.
\end{equation}

\begin{lemma}
\label{Lemma3-2}
There are $\e,C > 0$ such that if
$f : [0,\infty) \to \C$ satisfies 
\begin{equation}
\label{eq:small}
  | \lan f, h_1 \ran_X | \leq \e \| f \|_X
\end{equation}
then 
\begin{equation}
\| f \|_X^2 = 
  \int_0^{\infty}\big( |f_r|^2+\frac{m^2}{r^2} |f|^2 \big)rdr
  \leq C\int_0^{\infty} |L_0 f|^2 rdr.
\label{Lemma3-1:10}
\end{equation}
\end{lemma}
\begin{proof}
We may assume $f(r)$ is real-valued. Under our change of variables
$\tf(y) := f(e^{y/m})$ we have
\[
  \int_0^\infty (f_r^2 +\frac{m^2}{r^2}f^2)rdr = 
  m\int_\R ((\tf')^2 + \tf^2)  \, dy
\]
and
\[
  \int_0^{\infty} [L_0 f]^2 rdr
  = m \int_\R [\tilde{L}_0 \tilde{f}]^2  \, dy 
\]
where  
\[
  \tilde{L}_0 := \frac d{dy}+ \tanh y 
  = \sech(y) \frac{d}{dy} \frac{1}{\sech y}. 
\]
\newcommand{\snorm}[1]{\| #1 \|}
In the variable $y$, the {\it assumption} of the lemma becomes
\begin{equation}
\label{eq:approxorth}
  \lan \tilde{f}, \sech y \ran_{H^1(\R)} =
  \int_\R  \big({\tilde{f}}' (\sech y)'+ \tilde{f} \sech y\big)dy
  \leq \frac{\e}{\sqrt{m}} \snorm{\tf}_{H^1(\R)},
\end{equation}
and so it suffices to prove that~(\ref{eq:approxorth}) implies 
\begin{equation} 
\label{Lemma3-1:20}
  \| \tf \|_{H^1(\R)}^2 \le C \| \tilde{L}_0 \tilde{f} \|_{L^2(\R)}^2
  = C\lan \tf, \tilde{L}_0^* \tilde{L}_0 \tf \ran_{L^2(\R)}.
\end{equation}
The second-order differential operator 
\[
  H := \tilde{L}_0^* \tilde{L}_0 = - \frac{d^2}{dy^2} + 1 - 2\sech^2(y)
\]
is nonnegative with unique zero-eigenfunction (``ground state'') 
$\sech(y)$ (in fact, this operator is well-studied; see e.g. \cite{SS}).
Set $\phi := (1/\sqrt{2})\sech y$, so that $\|\phi\|_{L^2(\R)}=1$.
Write
\[ 
  \tilde{f} = a\phi + \varphi \quad \mbox{ with } \quad
  \lan \phi, \varphi \ran = 0.
\]
Similarly, decompose 
\[
  \phi = b(\phi''-\phi) + \psi \quad \mbox{ with } \quad
  \lan (\phi''-\phi), \psi \ran = 0. 
\]
Since 
$\phi''-\phi = -\sqrt{2}\sech^3 y$, $\psi$ is non-zero,
$b=-\|\phi^2\|_{L^2(\R)}^2 /(2\|\phi^3\|_{L^2(\R)}^2)$, 
and $\|\psi\|_{L^2} < 1$.
Now by~(\ref{eq:approxorth}),
\[
\begin{split}
  |a| &= | \lan \phi,\tilde{f} \ran| \leq 
  |\lan b(\phi''-\phi), \tilde{f} \ran| + |\lan \psi, \tilde{f} \ran| \\
  &= |b||\lan \phi, \tilde{f} \ran_{H^1}| + |\lan \psi, \tilde{f} \ran|
  \leq C\e \snorm{\tilde{f}}_{H^1(\R)}+| \lan \psi, \tilde{f} \ran|.
\end{split}
\] 
Using the above estimate and $(A+B)^2\leq q'A^2+qB^2$ where 
$q'$ is the H\"older conjugate of $q$ with $1<q<\infty$, we have
\[
\begin{split}
  \lan \tilde{f}, \tilde{f} \ran 
  &= |a|^2 + \| \varphi \|_{L^2}^2  
  \leq \big(C\e \snorm{\tilde{f}}_{H^1(\R)}+|\lan \psi, \tilde{f} \ran|\big)^2
  +\| \varphi \|_{L^2}^2 \\
  &\leq Cq'\e^2\snorm{\tilde{f}}^2_{H^1(\R)} + q |\lan \psi,\tilde{f} \ran|^2
  +\| \varphi \|_{L^2}^2 \\
  &\leq Cq'\e^2 \snorm{\tilde{f}}^2_{H^1(\R)}
  +q\| \psi \|_{L^2}^2 \|\tilde{f}\|_{L^2}^2 + \| \varphi \|_{L^2}^2.
\end{split}
\]
Choose $q$ to be such that $q \| \psi \|_{L^2}^2 <1$
to obtain
\[
  \| \tilde{f} \|_{L^2}^2 \leq C\e^2 \snorm{\tilde{f}}_{H^1(\R)} ^2+
  C\lan \tilde{f}, H \tilde{f} \ran,
\]
where we used 
$\| \varphi \|_{L^2}^2 \leq C \lan \varphi, H \varphi \ran
=C \lan \tilde{f}, H \tilde{f} \ran$.
On the other hand, we have 
\[
  \lan \tilde{f}',\tilde{f}' \ran \leq 
  \lan \tilde{f}, H\tilde{f}) + C\|\tilde{f}\|_{L^2}^2.
\]
Combining the two estimates above, we get
\[
  \| \tilde{f} \|_{H^1(\R)}^2 
  \leq C \lan \tilde{f}, H\tilde{f} \ran
  =C\int [\tilde{L}_0 \tilde{f}]^2 dy,
\]
provided $\e$ is sufficiently small.
Transforming back to the variable $r$, we obtain the estimate 
\eqref{Lemma3-1:10}, completing the proof.
\end{proof}

\medskip

{\bf Step 7: almost orthogonality.}
To apply the previous lemma to $z(r)$, we need to 
verify condition~(\ref{eq:small}), for which
we use the following lemma.

\begin{lemma}
\label{lem:orth}
For $z(r)$ and $\gamma(r)$ defined by \eqref{xi-def} and \eqref{xi-decomp},
\begin{equation}
\label{eq:orth1}
  \lan z_1, h_1 \ran_X = 
  \int_0^{\infty}(z_{1r} h_{1r} +\frac{m^2}{r^2}z_1h_1)rdr=0,
\end{equation}
\begin{equation}
\label{eq:orth2}
  \lan z_2, h_1 \ran_X = 
  \int_0^{\infty}(z_{2r} h_{1r} + \frac{m^2}{r^2} z_2h_1)rdr 
  = \int_0^{\infty} \frac{4m^2}{r^2}h^2_1h_3\gamma rdr.
\end{equation}
\end{lemma}
\begin{proof}
The pair $(\al(u),s(u))$ is the minimizer of
the differentiable function
\[
  F(s,\alpha) = \int_{\R^2}\abs{\nabla\big(e^{m\theta R}
  (v-h^{\alpha,s})\big)}^2dx.
\]
The lemma follows from the equations $\nabla_{\alpha,s} F(\al(u),s(u))
= 0$.
\end{proof}

\medskip

{\bf Step 8: proof of (b)}.
We will use the following abbreviation: 
\[
  ze := z_1 e + z_2 J^he;
  \quad\quad   z=z_1+iz_2.
\]
By~(\ref{eq:global}), we may choose $\delta$ sufficiently 
small so that  
\[
  \| u - e^{m\theta R} h^{s,\al} \|_{\dot H^1} < \delta_0
\] 
for any given $\delta_0 > 0$ (which will be specified later). 
So we have
\begin{equation}
\label{eq:xibound}
  \norm{\xi_r}^2_{L^2(\R^2)} + \norm{\frac{m}{r} R\xi}^2_{L^2(\R^2)} =
  \| u-e^{m\theta R}h^{s,\al} \|_{\dot H^1}^2
  < \delta^2_0.
\end{equation}

It is proved in Lemma~\ref{lem:pointwise} that~(\ref{eq:xibound})
implies the $L^\infty$ smallness of $\xi(r)$
for $\delta_0$ sufficiently small:
$\| \xi \|_{L^\infty} \leq C\delta_0$.
This immediately implies
\begin{equation}
\label{eq:smallness}
  \| z \|_{\infty} + \|\gamma \|_\infty \leq C\delta_0.
\end{equation}
Since $1=|v|^2 = |z|^2 + (1+\gamma)^2$ and $|\gamma|$ is small, 
\begin{equation}
\label{gamma-est}
\gamma = \sqrt{ 1 - |z|^2} -1\le 0, \quad 
|\gamma(r)| \le C |z(r)|^2, \quad
|\gamma_r(r)| \le C |z(r)z_r(r)|.
\end{equation}

By Lemma~\ref{lem:orth}, Eqns. (\ref{eq:smallness})--\eqref{gamma-est},
and $|h(r)| \le 1$,
\[
  \abs{\lan z, h_1 \ran_X} = 
  \abs{\int_0^{\infty}\frac{4m^2}{r^2}h^2_1h_3\gamma rdr} 
  \leq C\norm{z}_{\infty}\norm{\frac{h_1}{r}}_2\norm{\frac{z}{r}}_2
  \leq C\delta_0\norm{z}_X.
\]

Taking $\delta_0$ small enough so that $C\delta_0$
is less than $\e$ in Lemma~\ref{Lemma3-2}, we have
\begin{equation}
\label{eq:star}
\begin{split}
  \dist^2(u,\sO_m)
  &=2\pi\int_0^{\infty}\big(\abs{\xi_r}^2+\frac{m^2}{r^2}\abs{R \xi}^2 \big)rdr
  \leq C\int_0^{\infty}\big(\abs{z_r}^2+\frac{\abs{z}^2}{r^2}\big)rdr \\
  &\leq C\int_0^{\infty}\abs{L_0z}^2rdr.
\end{split}
\end{equation}

On the other hand, by \eqref{Eu-expand},
\[
  \delta_1^2 = \sE (u)-4\pi m = 
 \frac 12 \int _{\R^2}\abs{L\xi+\frac{m}{r}\xi_3\xi}^2 dx
  = \frac 12 \int _{\R^2} \abs{L(ze)+L(\gamma h)+\frac{m}{r}\xi_3\xi}^2 dx.
\]
Hence
\[
  \int _{\R^2} \abs{L_0z}^2dx=\int _{\R^2} \abs{L(ze)}^2dx
  \leq 3\delta_1^2
  + C\int_{\R^2} \big(\abs{L(\gamma h)}^2+\abs{\frac{m}{r}\xi_3\xi}^2\big)dx. 
\]
Using~(\ref{eq:smallness}) we have, for $\delta_0$ sufficiently small, 
\[
  \int_{\R^2}\abs{\frac{m}{r}\xi_3\xi}^2dx
  \leq C\norm{z}^2_{\infty}\int _{\R^2}\abs{\frac{z}{r}}^2dx.
\]
For the term $\|L(\gamma h)\|^2_{L^2}$, using \eqref{gamma-est}
we find
\[
  \int _{\R^2} |L(\gamma h)|^2\leq 
  C\int_{\R^2}\big(|\gamma_r h|^2+|\frac{\gamma h}{r}|^2
  +|\frac{\gamma}{r}h_3h|^2\big)
  \leq C\norm{z}_{\infty}^2\int_{\R^2}|L_0z|^2.
\]
Thus we get $(1-C\delta_0^2)\|L_0z\|^2_{L^2} \leq 4 \delta_1^2$.
Now choose $\delta_0 > 0$ sufficiently small (by choosing
$\delta$ small) so that $\|L_0 z\|^2_{L^2} \leq 5\delta_1^2$,
and therefore, $\dist(u,\sO_m) < C\delta_1$.
This completes the proof of 
part (b) of Theorem~\ref{thm:maps}.

\medskip

{\bf Step 9: proof of (c)}
By \eqref{Prop3-5:20} with $s=s(u)$ and $\al=\al(u)$, we have
$\int_{\R^2}\nabla u\cdot\nabla(e^{m\theta R}h^{s,\al})dx > 2\pi m$, 
since $\delta_1$ is small.
On the other hand, 
inequality~\eqref{Cor2-2:10} from Lemma~\ref{Lemma2-1} 
in Section~\ref{sec:lemmas} gives an upper bound: for any $\sigma \in (0,1)$,
\[
\begin{split}
  \int_{\R^2}\nabla u\cdot\nabla(e^{m\theta R}h^{s,\al})dx
  &\leq (\int_{\R^2}\frac{\abs{\nabla u}^2}{\abs{x}^{2\si}}dx)^{\frac{1}{2}}
  (\int_{\R^2}\abs{x}^{2\si}
  \abs{\nabla(e^{m\theta R}h^{s,\al})}^2dx)^{\frac{1}{2}} \\
  &\leq C\norm{\nabla u}^{1-\si}_{L^2(\R^2)}\norm{\nabla^2 u}^{\si}_{L^2(\R^2)}
  \norm{|x|^{\sigma}|\nabla (e^{m\theta R}h^{s,\al})| }_{L^2(\R^2)} \\
  &\leq C\norm{\nabla^2 u}^{\si}_{L^2(\R^2)}s^{\si}.
\end{split}
\]
Thus we obtain $\norm{\nabla^2 u}_{L^2(\R^2)}\geq C/s$,
completing the proof of Theorem~\ref{thm:maps}.
\end{proof}


\section{Global well-posedness vs. blow up}
\label{sec:flow}

In this section we complete the proof of our main result,
Theorem~\ref{mainthm:20}.

Let $u \in C([0,T);\dot H^2 \cap \Sigma_m)$ be a solution
of the Schr\"odinger flow equation~(\ref{int:100}) which 
conserves energy. We are assuming
\[
  \delta_1^2 := \sE(u_0) - 4\pi m = \sE(u(t)) - 4\pi m
  < \delta^2
\]
where $\delta$ is to be taken sufficiently small.
In particular, we choose $\delta$ small enough so
that Theorem~\ref{thm:maps} applies for each $t \in [0,T)$,
and so furnishes us with continuous functions
$s(t) \in (0,\infty)$ and $\alpha(t) \in \R$ such that
\begin{equation}
\label{eq:perturbation}
  \| u(x,t) - e^{(m \theta + \alpha(t)) R} h(r/s(t)) \|_{\dot H^1}
  < C \delta_1
\end{equation}
and
\[
  s(t) \| u(t) \|_{\dot H^2} > C > 0.
\]
Thus the first part of Theorem~\ref{mainthm:20}
is proved. It remains to show that 
\begin{equation}
\label{eq:goal}
  T < \infty \mbox{ and } 
  \lim_{t \to T^{-}} \| u(t) \|_{\dot H^2} = \infty \quad
  \implies \quad
  \liminf_{t \to T^-} s(t) = 0.
\end{equation}

Recall that we are writing
\[
  u(x,t) = e^{m \theta R} v(r), \quad\quad
  v(r) = e^{\al(t) R} [ h(r/s(t)) + \xi(r/s(t))],
\]
and~(\ref{eq:perturbation}) is equivalent to
\[
  \| e^{m \theta R} \xi(r) \|_{\dot H^1}^2
  = \| \xi_r(r) \|_{L^2}^2 + m^2 \norm{ \frac{R \xi(r)}{r} }_{L^2}^2
  < C \delta_1^2.
\]

To prove~(\ref{eq:goal}), we need estimates showing that 
$\|u(t)\|_{\dot H^2}$ is controlled as long as
$s(t)$ is bounded away from zero. 
These $\dot H^2$-estimates are obtained 
using the fact that the coordinates of
the tangent vector field
$v_r - \frac{m}{r}J^v Rv$,
with respect to a certain orthonormal frame,
satisfy a nonlinear Schr\"odinger-type equation,
and can be estimated using Strichartz estimates.

This construction, which was introduced in~\cite{CSU},
begins with a unit tangent vector field,
$\hat{e}(r) \in T_{v(r)} \S^2$, satisfying the
parallel transport condition 
\[
  D_r^v \hat{e} \equiv 0.
\]
Recall that $D_r^v$ is the covariant derivative,
acting on vector fields $\eta(r) \in T_{v(r)} \S^2$:
\[
  D_r^v \eta := P^v \eta_r
  = \eta_r - (v \cdot \eta_r)v
  = \eta_r + (v_r \cdot \eta)v.
\]
So $\hat{e}(r)$ and $J^v \hat{e}(r)$ form
an orthonormal frame on $T_v \S^2$.
Then $q(r) = q_1(r) + iq_2(r)$ is defined 
to be the coordinates of
$v_r - \frac{m}{r}J^v Rv \in T_v \S^2$ 
in this basis:
\[
  v_r - \frac{m}{r}J^v Rv    
  = q_1 \hat{e} + q_2 J^v \hat{e}.
\]
We will sometimes write 
$q\hat{e} := q_1\hat{e}+q_2J^v\hat{e}$ for convenience.
Note that by~(\ref{eq:1-11}),
\[
  \| q \|_{L^2(rdr)}^2 = \|  v_r - \frac{m}{r}J^v Rv \|_{L^2(rdr)}^2
  = \frac{1}{\pi}(\sE(u) - 4\pi m) = \delta_1^2/\pi
\]
is constant in time, and can be taken small. 

Define ${\nu}={\nu}_1+i{\nu}_2$ as follows:
\[
  J^v Rv = {\nu}_1\hat{e}+{\nu}_2J^v\hat{e}.
\] 
Again, we will sometimes denote 
${\nu}\hat{e} := {\nu}_1\hat{e}+{\nu}_2J^v\hat{e}$.
It is now straightforward, if somewhat involved,
to show that if $u(x,t)$ solves the Schr\"odinger map
equation~(\ref{int:100}), the complex function $q(r,t)$ 
solves the following nonlinear Schr\"odinger equation
with non-local nonlinearity 
(see \cite{CSU} for more details):
\begin{equation}
  iq_t=-\Delta_r q+\frac{(1-mv_3)^2}{r^2}q + \frac{m(v_3)_r}{r}q
  + qN(q),
\label{eq-4.8}
\end{equation}
where
\[
  N(q) ={\rm{Re}}\int_{r}^{\infty}(\bar{q}+\frac{m}{r}\bar{{\nu}})
  (q_r+\frac{1-mv_3}{r}q)dr.
\]
By changing variables to  $\tilde{q}:=e^{i(m+1)\theta}q$,
we obtain 
\begin{equation}
  i\tilde{q}_t+\Delta\tilde{q}-\frac{m(1+v_3)(mv_3-m-2)}{r^2}\tilde{q}
  -\frac{mv_{3,r}}{r}\tilde{q} - \tilde{q}N(q)=0. 
\label{eq-4.9}
\end{equation}
We will use this equation to obtain $H^1$ estimates on $q$.

For these estimates for $q$ to be useful,
we need to bound the original map $u(x,t)$
-- or equivalently $v(r,t)$ or $\xi(r,t)$ or $z(r,t)$ 
-- by $q$. Since $v(r,t) = e^{\al R}[h(r/s) + \xi(r/s,t)]$,
we have
\begin{equation}
  q\hat{e} = \frac{1}{s}
  e^{\al R}[L\xi + \frac{m}{r}\xi_3\xi]\left( \frac{r}{s} \right),
\label{Cor4-4:20}
\end{equation}
where, recall, 
\[ 
  L\xi = \xi_r + \frac{m}{r}(\xi_3 h+h_3 \xi),
  \quad \mbox{ and } \quad
  \xi = z_1 e + z_2 J^he + \gamma h.
\]
Since $\xi$ is small, we have, very roughly speaking,
\[
  \xi \approx ze, \quad \mbox{ and }
  L\xi \approx (L_0 z)e,
\]
and $z$ can be controlled by $L_0 z$.
More precisely, Lemma~\ref{Lemma4-5}, proved in the next section, 
gives the following bounds: for $2 \leq p < \infty$, and
provided $\delta$ is sufficiently small,
\begin{equation}
\label{eq:bound1}
  \norm{z_r}_{L^p(\R^2)} + \norm{\frac{z}{r}}_{L^p(\R^2)} 
  \leq C \left( s^{1-2/p} \norm{q}_{L^p(\R^2)} + 
  \norm{q}_{L^2(\R^2)} \right),
\end{equation}
\begin{equation}
\label{eq:bound2}
  \norm{z_{rr}}_{L^2(\R^2)} 
  \leq C \big( s\norm{q_r}_{L^2(\R^2)} + s\norm{\frac{q}{r}}_{L^2(\R^2)} 
  + s\norm{q}^2_{L^4(\R^2)} + \norm{q}_{L^2(\R^2)} \big),
\end{equation}
\begin{equation}
\label{eq:bound3}
  \| u \|_{\dot H^2(\R^2)} \leq C \big(
  \frac{1}{s} + \norm{q_r}_{L^2(\R^2)} + \norm{\frac{q}{r}}_{L^2(\R^2)} 
  + \norm{q}^2_{L^4(\R^2)} + \frac{1}{s} \norm{q}_{L^2(\R^2)} \big).
\end{equation}

We will use the following notation to denote space-time Lebesgue norms:
for an interval $I \subset \R$,
\[
  \| f \|^r_{L^r_t L^p_x(\R^2 \times I)} :=
  \int_I \left(\int_{\R^2} |f(x,t)|^p dx \right)^{r/p} dt.
\]  
We use~(\ref{eq:bound1})--(\ref{eq:bound2})
together with Strichartz estimates
and Equation~(\ref{eq-4.9}) to prove
the following estimates for $q$:
\begin{lemma}
\label{Lemma4-7}
For $\tau \geq 0$ and $\sigma > 0$, set
$I := (\tau, \tau + \sigma)$, $Q := \R^2 \times I$, and 
$X(Q) := L^4_t L^4_x(Q) \cap L^\infty_t L^2_x(Q)
\cap L^{8/3}_t L^8_x(Q)$.
Define $\underbar{s} := \inf_{t \in I} s(t)$.
If $\delta$ is sufficiently small, we have
\begin{equation}
\label{eq:q1}
  \norm{q}_{X(Q)} \leq 
  C \left( \norm{q_0}_{L^{2}_{x}(\R^2)}+
  (\underbar{s}^{-1}\sigma^{\frac{1}{2}}+\norm{q}^2_{L^4_x L^4_t(Q)})
  \norm{q}_{L^4_x L^4_t(Q)} \right),
\end{equation}
\begin{equation}
\label{eq:q2}
  \norm{\nabla\tilde{q}}_{X(Q)}
  \leq C( \norm{\nabla \tilde{q}(\tau)}_{L^{2}_{x}}
  +\big(
\underbar{s}^{-1}\sigma^{\frac{1}{2}} + 
\underbar{s}^{- \frac 32}\sigma^{\frac{3}{4}}
  +\norm{q}^2_{X(Q)}\big)
  \norm{\nabla\tilde{q}}_{L^\infty_tL^2_x(Q) \cap L^4_t L^4_x(Q)}).
\end{equation}
\end{lemma}

Before proving Lemma~\ref{Lemma4-7}, we show how it completes
the proof of our main theorem.

\medskip

\noindent
{\bf Completion of the proof of Theorem~\ref{mainthm:20}.}
We need to prove~(\ref{eq:goal}), so suppose that 
\begin{equation}
\label{eq:blowup}
  \lim_{t \to T^-} \| u(t) \|_{\dot H^2} = \infty
\end{equation} 
and
\begin{equation}
\label{eq:s(t)}  
  \liminf_{t\rightarrow T^{-}} s(t) = s_0 > 0.
\end{equation}
Our goal is to derive a contradiction. 
By~(\ref{eq:s(t)}), we have 
$s(t) \geq s^*$ for all $0 \leq t < T$, for some $s^* > 0$.
So we may take $\underbar{s} = s^*$ in the estimates~(\ref{eq:q1}-\ref{eq:q2})
for any time interval $I \subset [0,T)$.
If $\sigma$ is sufficiently small depending on $s^*$
($\sigma^{1/2}  < s^*/2C$), and 
$\| q_0 \|_{L^2} \leq C \delta_1$ is taken sufficiently small,
estimate~(\ref{eq:q1}) implies
\[
  \| q \|_{X(Q)} \leq C \delta_1.
\]
Using this estimate in~(\ref{eq:q2}), for $\delta_1$
and $\sigma$ sufficiently small we obtain
\[
  \| \nabla \tilde{q} \|_{L^\infty_t L^2_x(Q)}
  \leq C \| \nabla \tilde{q}(\tau) \|_{L^2}.
\]
In particular, taking $\tau$ close to $T$, we see
\[
  \limsup_{t \to T^-} \| \nabla \tilde{q}(t) \|_{L^2} < \infty.
\]
Then using~(\ref{eq:bound3}),
and $|q_r| + |q/r| \leq C |\nabla \tilde{q}|$,
we find $\limsup_{t \to T^-} \| u(t) \|_{\dot H^2} < \infty$,
contradicting~(\ref{eq:blowup}). 
This completes the proof of Theorem~\ref{mainthm:20}.
$\Box$

\medskip

\noindent
{\bf Proof of Lemma~\ref{Lemma4-7}.}
Strichartz estimates for the inhomogeneous Schr\"odinger equation
(see eg. \cite{C}), applied to~(\ref{eq-4.9}), give
\begin{equation}
\label{eq:strich}
  \| q \|_{X(Q)} \leq C( \| q(\tau) \|_{L^2}
  + \| F \|_{L^{4/3}_t L^{4/3}_x (Q)} )
\end{equation}
where
\[
  F := \frac{m(1+v_3)(mv_3-m-2)}{r^2}\tilde{q}
  +\frac{m (v_3)_r}{r}\tilde{q}+\tilde{q}N(q).  
\]
Conservation of the $L^2$-norm of $q(t)$
(equivalent to conservation of energy $\sE(u)$)
means we can replace $\| q(\tau) \|_{L^2}$
by $\| q_0 \|_{L^2}$ in~(\ref{eq:strich}).
Using $v_3(r) = h_3(r/s) +\xi_3(r/s)$ and
$\xi_3=z_2h_1+\gamma h_3$, we find
\begin{equation}
\begin{split}
  \norm{\frac{1+v_3}{r^2}}_{L^2_x} 
  &= \frac{1}{s} \norm{\frac{1+h_3+z_2h_1+\gamma h_3}{r^2}}_{L^2_x}
  \leq \frac{C}{s} \left( \norm{\frac{1+h_3}{r^2}}_{L^2_x}
  +\norm{\frac{z_2h_1}{r^2}}_{L^2_x}+\norm{\frac{\gamma h_3}{r^2}}_{L^2_x}\right) \\
  &\leq \frac{C}{s}(1 + \norm{\frac{z}{r}}_{L^2_x}
  +\norm{\frac{z}{r}}^2_{L^4_x})
  \leq C(s^{-1}(1 + \norm{q}_{L^2} + \norm{q}^2_{L^2})
  +\norm{q}^2_{L^4}).
\end{split}
\label{Lemma4-7:10}
\end{equation}
Using~\eqref{Lemma4-7:10} and the uniform 
boundedness of $\|q\|_{L^2}$, we get
\[
\begin{split}
  \norm{\frac{m(1+v_3)(mv_3-m-2)}{r^2}\tilde{q}}_{L^{4/3}_t L^{4/3}_x}
  &\leq C\norm{\frac{1+v_3}{r^2}\tilde{q}}_{L^{4/3}_t L^{4/3}_x}
  \leq C\norm{ \norm{\frac{1+v_3}{r^2}}_{L^2_x}
  \norm{q}_{L^4_x}}_{L^{4/3}_t} \\
  &\leq C(\underbar{s}^{-1}{\sigma}^{\frac{1}{2}}+\norm{q}^2_{L^4_t L^4_x})
  \norm{q}_{L^4_t L^4_x}.
\end{split}
\]
Next we estimate $\|(m(v_3)_r/r)q\|_{L^{4/3}_t L^{4/3}_x}$.
Compute
\[
  \frac{1}{r} (v_3)_r = \frac{1}{sr}
  \left[ \frac{m h^2_1}{r}+ h_1 (z_2)_r 
  -\frac{mh_1h_3z_2}{r}+\gamma_r h_3+\frac{m\gamma h^2_1}{r}
  \right]\left( \frac{r}{s} \right).
\]
({\it Remark:} the notation means evaluate the quantity
in the square brackets with $r$ replaved by $r/s$.
We use this notation frequently in what follows.) 
Therefore, using boundedness of 
$z$, $h_3$, and $h_1/r$, and
$\gamma=O(|z|^2)$, $\gamma_r=O(|zz_r|)$,
we have
\begin{equation}
  \norm{\frac{(v_3)_r}{r}}_{L^2}
  \leq \frac{C}{s} \big( \norm{\frac{h^2_1}{r^2}}_{L^2}
  + \norm{z_{r}}_{L^2} + \norm{\frac{z}{r}}_{L^2} \big) 
  \leq Cs^{-1}(1+\norm{q}_{L^2}).
\label{Lemma4-7:20}
\end{equation}
So
\[
  \norm{\frac{m(v_3)_r}{r}q}_{L^{4/3}_t L^{4/3}_x}\leq 
  C\norm{\norm{\frac{v_{3,r}}{r}}_{L^2_x}\norm{q}_{L^4_x}}_{L^{4/3}_t}
  \leq C(1+\norm{q_0}_{L^2})\underbar{s}^{-1}{\sigma}^{\frac{1}{2}}
  \norm{q}_{L^{4}_t,L^{4}_x}.
\]
Next we need to estimate the nonlocal term
\[
\begin{split}
  N(q) &= {\rm{Re}}\int^{\infty}_r(\bar{q}+\frac{m\bar{{\nu}}}{r})
  (q_r+\frac{1-mv_3}{r}q)dr \\
  &= -\frac{|q|^2}{2} + {\rm{Re}}\int^{\infty}_r
  \big(\frac{1-mv_3}{r}|q|^2+\frac{m\bar{{\nu}}}{r}q_r
  +\frac{m\bar{{\nu}}(1-mv_3)}{r^2}q\big)dr \\
  &:= -\frac{|q|^2}{2} + (R_1+R_2+R_3).
\end{split}
\]
First note
\[
  \| |q|^2q \|_{L^{4/3}_t L^{4/3}_x}
  = \| q \|_{L^4_t L^4_x}^3.
\]
Now consider $\tilde{q}R_1$.
Using the estimate 
\[
  \|f(r)\|_{L^2} \leq C \|r f_r(r)\|_{L^2}
\]
(Hardy's inequality in $\R^4$ for radial functions), we have
\[
  \norm{\tilde{q}R_1}_{L^{\frac{4}{3}}_t,L^{\frac{4}{3}}_x}
  \leq \|q\|_{L^4L^4} \|R_1\|_{L^2L^2}
  \leq C\norm{q}^3_{L^4_t L^4_x}.
\]
Next we note that $|\bar{{\nu}}|=|J^vRv|=\sqrt{1-v^2_3}=|Rv|$. 
We consider next $\tilde{q}R_3$. By Hardy again,
\begin{equation}
\begin{split}
  \norm{\tilde{q}R_3}_{L^{\frac{4}{3}}_x}
  &\leq C\norm{\tilde{q}\int^{\infty}_r\frac{m\bar{{\nu}}(1-mv_3)}{r^2}
  qdr}_{L^{\frac{4}{3}}}
  \leq C\norm{\frac{\bar{{\nu}}}{r}q}_{L^2}\norm{q}_{L^4}
  \leq C\norm{\frac{\bar{{\nu}}}{r}}_{L^4}\|q\|^2_{L^4} \\
  &\leq C\norm{\frac{1-v^2_3}{r^2}}^{\frac{1}{2}}_{L^2}\norm{q}^2_{L^4}
  \leq C\norm{\frac{1+v_3}{r^2}}^{\frac{1}{2}}_{L^2}\norm{q}^2_{L^4}
  \leq C(s^{-\frac{1}{2}}+\norm{q}_{L^4})\norm{q}^2_{L^4},
\label{Lemma4-7:30}
\end{split}
\end{equation}
where we used~\eqref{Lemma4-7:10} again.
Thus
\begin{equation}
  \norm{\tilde{q}R_3}_{L^{\frac{4}{3}}_t,L^{\frac{4}{3}}_x}
  \leq C\norm{ (s^{-\frac{1}{2}}+\norm{q}_{L^4_x})
  \norm{q}^2_{L^4_x}}_{L^{\frac{4}{3}}_t}
  \leq C(\underbar{s}^{-\frac{1}{2}}{\sigma}^{\frac{1}{4}}
  +\norm{q}_{L^4_t L^4_x})\norm{q}^2_{L^4_t L^4_x}.
\label{Lemma4-7:40}
\end{equation}
It remains to estimate $\tilde{q}R_2$.
We rewrite $R_2$, using integration by parts:
\[
\begin{split}
  R_2 &= - {\rm{Re}}\frac{m}{r}\bar{{\nu}}q
  -{\rm{Re}}\int^{\infty}_r\left(\frac{m}{r}\bar{{\nu}}_r
  -\frac{m}{r^2}\bar{{\nu}}\right)qdr \\
  &=- {\rm{Re}}\frac{m}{r}\bar{{\nu}}q
  +{\rm{Re}}\int^{\infty}_r\left(\frac{m v_3}{r}\bar{q}
  +\frac{m(1+mv_3)}{r^2}\bar{{\nu}}\right)q dr,
\end{split}
\]
where we used $\bar{{\nu}}_r=-v_3\bar{q}-(m/r)v_3\bar{{\nu}}$.
So
\[
  \norm{\tilde{q}R_2}_{L^{\frac{4}{3}}_t,L^{\frac{4}{3}}_x}
  \leq C \left( \norm{q\frac{\bar{{\nu}}}{r}\bar{q}}_{L^{\frac{4}{3}}_t,
  L^{\frac{4}{3}}_x}
  +\norm{q \int^{\infty}_r\frac{|q|^2}{r}
  dr}_{L^{\frac{4}{3}}_t,L^{\frac{4}{3}}_x}
  +\norm{q \int^{\infty}_r\frac{\bar{{\nu}}}{r^2}
  qdr}_{L^{\frac{4}{3}}_t,L^{\frac{4}{3}}_x} \right).
\]
The last term is estimated as follows:
\begin{equation}
\label{Lemma4-7:50}
\begin{split}
  \norm{q\int^{\infty}_r\frac{\bar{{\nu}}}{r^2}
  qdr}_{L^{\frac{4}{3}}_t,L^{\frac{4}{3}}_x}
  &\leq C\norm{\norm{\frac{\bar{{\nu}}}{r}q}_{L^2_x}\|q\|_{L^4_x}
  }_{L^{\frac{4}{3}}_t}
  \leq C\norm{ \norm{\frac{\bar{{\nu}}}{r}}_{L^4_x}\norm{q}^2_{L^4_x}
  }_{L^{\frac{4}{3}}_t} \\
  &\leq C\norm{(s^{-\frac{1}{2}}+\norm{q}_{L^4_x})\norm{q}^2_{L^4_x}
  }_{L^{\frac{4}{3}}_t} \\
  &\leq C(\underbar{s}^{-\frac{1}{2}}{\sigma}^{\frac{1}{4}}
  +\norm{q}_{L^4_t L^4_x})\norm{q}^2_{L^4_t L^4_x},
\end{split}
\end{equation}
where we used the same computations as in~\eqref{Lemma4-7:30} 
and~\eqref{Lemma4-7:40}. The first term 
can be treated in a similar manner, leading to the same
estimate as in~\eqref{Lemma4-7:50}. 
The estimate for the second term has been done already.

Returning now to~(\ref{eq:strich}), and using
the above estimates, we 
have establishesed~(\ref{eq:q1}).

Next we need to estimate the derivative of $\tilde{q}$
in order to establish~(\ref{eq:q2}).
Denote $w := \partial_{x_i} \tilde{q}$ for $i=1,2$. 
Then $w$ solves 
\begin{equation}
\label{eq:w}
\begin{split}
  iw_t + \Delta w &= \frac{m(1+v_3)(mv_3-m-2)}{r^2}w +
  \left( \frac{m(1+v_3)(mv_3-m-2)}{r^2} \right)_{x_i} \tilde{q} \\
  & \quad + \frac{m(v_3)_r}{r}w + \left( \frac{m(v_3)_r}{r} \right)_{x_i}
  \tilde{q} + N(q)w + N(q)_{x_i}\tilde{q}.
\end{split}
\end{equation}
Using the previous estimates, we can estimate the various terms 
involving $w$ in the right hand side:
\[
\begin{split}
  \noindent \bullet\quad 
  \norm{\frac{m(1+v_3)(mv_3-m-2)}{r^2}w}_{L^{4/3}_t L^{4/3}_x}
  &\leq C\norm{(s^{-1}+\norm{q}^2_{L^4_x})\norm{w}_{L^4_x}}_{L^{4/3}_t} \\
  &\leq C( \underbar{s}^{-1}{\sigma}^{\frac{1}{2}}+\norm{q}^2_{L^4_t L^4_x})
  \norm{\nabla\tilde{q}}_{L^4_t L^4_x}
\end{split}
\]
\[
  \noindent \bullet\quad \norm{\frac{m(v_3)_r}{r}w}_{L^{4/3}_t L^{4/3}_x}
  \leq C\norm{s^{-1}\norm{w}_{L^4_x}}_{L^{4/3}_t}
  \leq C \underbar{s}^{-1}{\sigma}^{\frac{1}{2}}
  \norm{\nabla\tilde{q}}_{L^4_t L^4_x}
\]
\[
  \noindent \bullet\quad \norm{N(q)w}_{L^{4/3}_t L^{4/3}_x}
  \leq C(\underbar{s}^{-\frac{1}{2}}{\sigma}^{\frac{1}{4}}+
  \norm{q}_{L^4_t L^4_x}) \norm{q}_{L^4_t L^4_x}
  \norm{\nabla\tilde{q}}_{L^4_t L^4_x}
\]
Now the other terms.
First note that 
$|\nabla \tilde{q}|^2 \sim |\tilde{q}|^2+|\tilde{q}/r|^2$,
and thus $\|\tilde{q}/r\|_{L^p(\R^2)}\leq 
C\|\nabla \tilde{q}\|_{L^p(\R^2)}$ for any $1\leq p\leq\infty$.
Due to \eqref{Lemma4-7:10} and \eqref{Lemma4-7:20}, we have
\[
\begin{split}
  \norm{ \left( \frac{(1+v_3)(mv_3-m-2)}{r^2} \right)_{x_i}
  \tilde{q}}_{L^{4/3}_t L^{4/3}_x}
  &\leq C \left( \norm{\frac{(v_3)_r}{r^2}
  \tilde{q}}_{L^{4/3}_t L^{4/3}_x}
  +\norm{\frac{(1+v_3)}{r^3}\tilde{q}}_{L^{4/3}_t L^{4/3}_x} \right) \\
  &\leq C\norm{ \left( \norm{\frac{(v_3)_r}{r}}_{L^2_x}
  +\norm{\frac{1+v_3}{r^2}}_{L^2_x} \right)
  \norm{\frac{\tilde{q}}{r}}_{L^4_x}}_{L^{\frac{4}{3}}_t} \\
  &\leq C\norm{(s^{-1}+\norm{q}^2_{L^4_x})
  \norm{\nabla\tilde{q}}_{L^4_x}}_{L^{\frac{4}{3}}_t} \\
  &\leq C(\underbar{s}^{-1}{\sigma}^{\frac{1}{2}}+\norm{q}^2_{L^4_t L^4_x})
  \norm{\nabla\tilde{q}}_{L^4_t L^4_x}.
\end{split}
\]
Next we consider 
\[
  \norm{ \left( \frac{m(v_3)_r}{r} \right)_{x_i} 
  \tilde{q}}_{L^{4/3}_t L^{4/3}_x}
  \leq C \left( \norm{\frac{(v_3)_r}{r^2}
  \tilde{q}}_{L^{\frac{4}{3}}_t L^{\frac{4}{3}}_x}
  +\norm{\frac{(v_3)_{rr}}{r}\tilde{q}}_{L^{4/3}_t L^{4/3}_x} \right).
\]
The first term in the right side can be estimated as before:
\[
  \|((v_3)_r/r^2)\tilde{q}\|_{L^{4/3}_t L^{4/3}_x}
  \leq C \underbar{s}^{-1}{\sigma}^{1/2}\|\nabla\tilde{q}\|_{L^4_t L^4_x}.
\] 
Recalling $(h_1)_r = -(m/r)h_1h_3$, 
$(h_3)_r = (m/r)h^2_1$, and $h_1/r, h_3$ bounded, we find
\[
\begin{split}
  |(v_3)_{rr}| &\leq \frac{1}{s^2}
  \left| m\left(\frac{h_1^2}{r}\right)_r + 
  (h_1z_{2,r})_r -m \left(\frac{h_1h_3z_2}{r} \right)_r
  +(\gamma_rh_3)_r + m \left( \frac{\gamma h^2_1}{r} \right)_r \right|(r/s) \\
  &\leq \frac{C}{s^2}\left[ |(\frac{h_1^2}{r})_r| +|\frac{h_1z_2}{r^2}| 
  +|\frac{h_1(z_2)_r}{r}|+|h_1z_{rr}|+|zz_{rr}|+|z_{r}|^2 \right](r/s).
\end{split}
\]
We estimate term by term:
\[
\noindent \bullet\quad 
\norm{\frac{1}{s^2} \left[ \left(\frac{h_1^2}{r}\right)_r 
\left(\frac{r}{s}\right) \right] 
\frac{q}{r}}_{L^{4/3}_t L^{4/3}_x}
\leq C\norm{\frac{1}{s^2} \left( \left|\frac{h_1^2}{r^2}
+ \frac{h^2_1|h_3|}{r^2} \right| \left(\frac{r}{s}\right) \right)
\frac{q}{r}}_{L^{4/3}_t L^{4/3}_x}
\leq C \underbar{s}^{-1}{\sigma}^{\frac{1}{2}}
\norm{\nabla\tilde{q}}_{L^4_t L^4_x}.
\] 
\[
\begin{split}
\noindent \bullet\quad 
\norm{\frac{1}{s^2} \left( \left| \frac{h_1|z_2|}{r^2} +
\frac{h_1|(z_2)_r|}{r} \right| \left(\frac{r}{s}\right)
\right) \frac{q}{r}}_{L^{4/3}_t L^{4/3}_x}
&\leq C\norm{ s^{-1}(\norm{\frac{z}{r}}_{L^2_x}+\norm{z_r}_{L^2_x})
\norm{\frac{q}{r}}_{L^4_x}}_{L^{\frac{4}{3}}_t} \\
&\leq C\underbar{s}^{-1}{\sigma}^{\frac{1}{2}}
\norm{\nabla\tilde{q}}_{L^4_t L^4_x}.
\end{split}
\]
\[
\noindent \bullet\quad 
\norm{\frac{1}{s^2} \left( |z_{r}|^2 \left(\frac{r}{s}\right) \right)
\frac{q}{r}}_{L^{4/3}_t L^{4/3}_x}
\leq C\norm{ \frac{1}{s}\norm{z_r}^2_{L^4_x}
\norm{\frac{q}{r}}_{L^4_x}}_{L^{\frac{4}{3}}_t}   
\leq C(\norm{q}^2_{L^4_t L^4_x} + \underbar{s}^{-1} \sigma^{1/2})  
\norm{\nabla\tilde{q}}_{L^4_t L^4_x}.
\]
For the  remaining term,  using \eqref{eq:bound2},
\begin{align*}
\bullet\quad 
\norm{\frac{1}{s^2} \left( [(|h_1|+|z|) z_{rr}]
\left(\frac{r}{s}\right) \right) 
\frac{q}{r}} 
&  _{L^{4/3}_t L^{4/3}_x}
\leq C\norm{ \frac{1}{s} \norm{z_{rr}}_{L^2_x} 
\norm{(|h_1|+|z|)\bke{\frac rs} \frac{q}{r}}_{L^4_x}}_{L^{\frac{4}{3}}_t} \\
&\leq C\norm{ (\norm{q_r}_{L^2_x}+\norm{\frac{q}{r}}_{L^2_x}) 
\norm{(|h_1|+|z|)\bke{\frac rs}\frac{q}{r}}_{L^4_x}}_{L^{\frac{4}{3}}_t} \\
& \quad +C\norm{ (\norm{q}^2_{L^4_x}+s^{-1}\norm{q}^2_{L^2_x}) 
\norm{(|h_1|+|z|)\bke{\frac rs}\frac{q}{r}}_{L^4_x}}_{L^{\frac{4}{3}}_t} \\
&\leq C(\norm{q}^2_{L^{\frac{8}{3}}_t L^8_x}+\underbar{s}^{-\frac{3}{2}}
{\sigma}^{\frac{3}{4}})\norm{\nabla\tilde{q}}_{L^{\infty}_t L^2_x}\\ 
&\quad +C(\norm{q}^2_{L^4_t L^4_x}+ \underbar{s}^{-1}
{\sigma}^{\frac{1}{2}}) \norm{\nabla \tilde{q}}_{L^4_t L^4_x},
\end{align*}
where we used, first by \eqref{eq:bound1} with $p=8$,
\[
\norm{(|h_1|+|z|)(r/s)\frac{q}{r}}_{L^4_x} 
\le C \norm{\frac{(|h_1|+|z|)(r/s)} r}_{L^8_x} \norm{q}_{L^8_x}
\le C \norm{q}_{L^8_x}^2 + C s^{-3/2} (1+\norm{q}_{L^2_x}^2),
\]
and then $\norm{(|h_1|+|z|)(r/s)\frac{q}{r}}_{L^4_x} \le C \norm{q/r}_{L^4_x}$.

It remains to estimate $N(q)_{x_i}q$:
\[
\begin{split}
\norm{N(q)_{x_i}\tilde{q}}_{L^{4/3}_t L^{4/3}_x}
&\leq C\norm{(|q|^2+|\frac{\bar{\nu}}{r}q|)
(|q_r|+|\frac{q}{r}|)}_{L^{4/3}_t L^{4/3}_x} \\
&\leq C\norm{ \norm{(|q|^2+|\frac{\bar{\nu}}{r}q|)}_{L^2_x}
\norm{(|q_r|+|\frac{q}{r}|)}_{L^4_x}}_{L^{\frac{4}{3}}_t} \\
&\leq C\norm{(\norm{q}^2_{L^4_x}
+\norm{\frac{\nu}{r}}^2_{L^4_x})\norm{\nabla\tilde{q}}_{L^4_x}
}_{L^{\frac{4}{3}}_t}
\leq C( \underbar{s}^{-1}{\sigma}^{\frac{1}{2}}
+\norm{q}^2_{L^4_t L^4_x})\norm{\nabla\tilde{q}}_{L^4_t L^4_x}.
\end{split}
\]
Now applying Strichartz estimates to~(\ref{eq:w}),
and using the estimates established above, we obtain
\[
\norm{\nabla\tilde{q}}_{L^{\infty,2}_{t,x}\cap L^4_t L^4_x}
\leq \norm{\nabla \tilde q(\tau)}_{L^{2}_{x}}
+C((\underbar{s}^{-1}+\underbar{s}^{-2}){\sigma}^{\frac{1}{2}}
+\| q \|_{X(Q)}^2)
\norm{\nabla\tilde{q}}_{L^{\infty,2}_{t,x}\cap L^4_t L^4_x}.
\]
This completes the proof of Lemma~\ref{Lemma4-7}.
$\Box$


\section{Technical lemmas}
\label{sec:lemmas}

In this section we collect some of the technical lemmas
used in the proof of the main theorem in the previous 
sections.

\subsection{Some inequalities for radial functions}

We begin with some inequalities for radial functions.

\begin{lemma} 
\label{Lemma2-1}
\begin{enumerate}
\item
Let $0\le \si <1$, and suppose $f\in H^1(\R^2)$ is radial. Then
\begin{equation}
  \int _0^\infty \frac {f^2}{r^{2\si}}\, r \, dr \le C_\si 
  \bke{\int_0^\infty f^2\, r \, dr}^{1-\si}
  \bke{\int _0^\infty |f_r|^2\, r\, dr}^{\si}.
\label{Lemma2-1:5}
\end{equation}
$\lim_{\si \to 1-} C_\si = \infty$ and the estimate is false
if $\si=1$.
\item
Let $0\leq\sigma<1$, and suppose $f\in H^1(\R^2)$. Then
\begin{equation}
  \int_{\R^2}\frac{f^2}{|x|^{2\si}}dx
  \leq C_{\sigma}\bke{\int_{\R^2} |f|^2 dx}^{1-\si}
  \bke{\int_{\R^2} |\nabla f|^2 dx}^{\si}.
\label{Cor2-2:10}
\end{equation}
\item
Suppose $f\in H^1_{\rm{loc}}(\R^2)$ is radial with 
$f_r, f/r \in L^2(\R^2)$. Then 
\begin{equation}
  \norm{f}_{L^{\infty}(\R^2)}\leq C\bke{\int_0^{\infty}(|f_r|^2
  +\frac{f^2}{r^2})rdr}.
\label{Lemma2-1:7}
\end{equation}
\end{enumerate}
\end{lemma}
\begin{proof}
We first show that 
\begin{equation}
\int_0^{\infty}\frac{f^2(r)}{r^{2\si}}rdr\leq 
\frac{1}{1-\si}\bke{\int_0^{\infty}\frac{f^2(r)}{r^{4\si-2}}rdr}^{\frac{1}{2}}
\bke{\int_0^{\infty}(f_r(r))^2rdr}^{\frac{1}{2}}.
\label{Lemma2-1:10}
\end{equation}
Indeed, by changing the order of integration, we get 
\[
\begin{split}
  \int_0^{\infty}\frac{f^2(r)}{r^{2\si}}rdr
  &=-\int_0^{\infty}\frac{r}{r^{2\si}}\int_r^{\infty}[f^2(s)]_sdsdr
  =-2\int_0^{\infty}f(s)f_s(s)ds\int_0^s\frac{r}{r^{2\si}}dr \\
  &=-\frac{1}{1-\si}\int_0^{\infty}s^{2-2\si}f(s)f_s(s) ds \\
  &\leq \frac{1}{1-\si}\bke{\int_0^{\infty}s^{2-4\si}f^2(s)sds}^{\frac{1}{2}}
\bke{\int_0^{\infty}(f_s(s))^2sds}^{\frac{1}{2}}.
\end{split}
\]
In particular, if $\si=1/2$, \eqref{Lemma2-1:10} immediately implies 
\eqref{Lemma2-1:5}. We note also that the estimate \eqref{Lemma2-1:5} is 
immediate in the case $\sigma=0$. 
Let $\si_i=1-(1/2^i)$ where $i \ge 0$ is an integer.
From the estimate \eqref{Lemma2-1:10} with $\si = \si_{i+1}$, we have
\[
  \int_0^{\infty}\frac{f^2(r)}{r^{2\si_{i+1}}}rdr\leq
  2^{i+1}\bke{\int_0^{\infty}\frac{f^2(r)}{r^{2\si_i}}rdr}^{\frac{1}{2}}
  \bke{\int_0^{\infty}(f_r(r))^2rdr}^{\frac{1}{2}}, \quad i=0,1,2,...\,\,. 
\]
Iterating this estimate, we obtain
\begin{equation}
  \int_0^{\infty}\frac{f^2(r)}{r^{2\si_{n+1}}}rdr\leq
  C_{n+1}\bke{\int_0^{\infty}f^2(r)rdr}^{\frac{1}{2^{n+1}}}
  \bke{\int_0^{\infty}(f_r(r))^2 rdr}^{1-\frac{1}{2^{n+1}}},
\label{Lemma2-1:20}
\end{equation}
for some constant $C_{n+1}$. (One can solve $C_{n+1} = 2^{n+1} \sqrt
{C_n}$ and $C_0=1$ to get $C_{n+1} = 2^{2n + 2^{-n}}$, which is
certainly not the best constant.)
It remains to consider the general case $\si \in [0,1)$.
Let $k\geq 0$ be an integer with $\si_k\leq
\si<\si_{k+1}$.  There exists $0\leq\theta<1$ such that
$\si=\theta\si_{k+1}+(1-\theta)\si_{k}=1-(2-\theta)/2^{k+1}$.  Using
the H\"older inequality and \eqref{Lemma2-1:20}, we get
\[
\begin{split}
\int_0^{\infty}\frac{f^2}{r^{2\si}}rdr &\leq 
\bke{\int_0^{\infty}\frac{f^2}{r^{2\si_{k+1}}}rdr}^{\theta}
\bke{\int_0^{\infty}\frac{f^2}{r^{2\si_k}}rdr}^{1-\theta} \\
&\leq C_{\si}\bke{\int_0^{\infty} f^2(r)rdr}^{1-\si}
\bke{\int_0^{\infty} (f')^2(r)rdr}^{\si},
\end{split}
\]
where $C_{\si}=C_{k+1}^{\theta}C_{k}^{1-\theta}$.
This completes the proof of the first estimate \eqref{Lemma2-1:5}. 

To see that this estimate fails at the endpoint $\sigma=1$,
fix a smooth, non-negative, non-decreasing function $\eta(r)$,
supported in $(1/2,\infty)$, and with $1-\eta(r)$ supported in
$[0,3/2)$. Then it easy to check that 
$f_{\delta}(r) := \eta(r/\delta) - \eta(r)$
provides a counterexample to the endpoint estimate as $\delta \to 0$.
Note that $f_{\delta}(0) = 0$ for all $\delta$. 

The second estimate~\eqref{Cor2-2:10} is an
immediate consequence of~\eqref{Lemma2-1:5}.
Using polar coordinates, we obtain
\[
  \int_0^{2\pi}\int_0^{\infty}
  \frac{|f(r,\theta)|^2}{r^{2\si}}r drd\theta
  \leq C_{\sigma}\int_0^{2\pi}
  \big(\int_0^{\infty} |f(r,\theta)|^2rdr\big)^{1-\si}
  \big(\int_0^{\infty} |\pd_r f(r,\theta)|^2rdr\big)^{\si}d\theta.
\]
Using $|\pd_r f|^2\leq |\nabla f|^2$, we obtain  estimate 
\eqref{Cor2-2:10} from H\"older's inequality.

For the third estimate \eqref{Lemma2-1:7},
we introduce the new variable $y$ defined by $r=e^y$ and 
denote $g(y)=f(e^y)=f(r)$. Then it is immediate that
\[
  \int_0^{\infty}\big(|f_r(r)|^2 +\frac{f^2(r)}{r^2}\big)rdr
  =\int_{-\infty}^{\infty}\big(|g'(y)|^2+|g(y)|^2\big)dy.
\]
By Sobolev embedding, we have 
$\| g \|_{L^{\infty}(\R)} \leq C \| g \|_{H^1(\R)}$. 
Transforming back to the original variable completes our proof.
\end{proof}


\begin{lemma}
\label{Lemma4-1}
Let $g:\R^2\rightarrow{\mathbb C}$ be radial and bounded with $g,g'
\in L^p_{loc}$, $2 < p <\infty$. Assume $(\partial_r-\frac{m}{r})g(r)
\in L^p(\R^2)$ for some $m \geq 1$. Then $g(r)/r \in L^p(\R^2)$ and
\[
  \norm{\frac{g(r)}{r}}_{L^p(\R^2)}\leq 
  C\norm{(\partial_r-\frac{m}{r})g(r)}_{L^p(\R^2)}.
\]
\end{lemma}
\begin{proof}
Let $0<r_1<r_2<\infty$ and denote $A= \{x\in\R^2:r_1<|x|<r_2\}$.
Consider
\[
  I:= - 2 \pi \Re \int_{r_1}^{r_2}(g_r-\frac{m}{r}g)|\frac{g}{r}|^{p-2}
  \frac{\bar{g}}{r}rdr.
\]
On one hand, $I\leq C\|g/r\|^{p-1}_{L^p(A)}\|g_r-\frac{m}{r}g\|_{L^p(A)}$
by H\"older inequality. On the other hand, 
\[
\begin{split}
  I &= m\norm{\frac{g}{r}}^p_{L^p(A)}
  -\frac{2\pi}{p}\frac{|g(r)|^p}{r^{p-2}}\big|^{r_2}_{r_1}
  -\frac{p-2}{p}\norm{\frac{g}{r}}^p_{L^p(A)} \\
  &\geq (m-1+\frac{2}{p})\norm{\frac{g}{r}}^p_{L^p(A)}
  -\frac{2\pi}{p}\frac{|g(r_2)|^p}{r_2^{p-2}}.
\end{split}
\]
Thus
\[
  (m-1+\frac{2}{p})\norm{\frac{g}{r}}^p_{L^p(A)}
  \leq \frac{2\pi}{p}\frac{|g(r_2)|^p}{r_2^{p-2}} +
  \norm{\frac{g}{r}}^{p-1}_{L^p(A)}
  \norm{g_r-\frac{m}{r}g}_{L^p(A)}.
\]
This gives a bound for $\norm{g(r)/r}_{L^p(A)}$ uniformly in
$r_1,r_2$. Hence $g(r)/r \in L^p(\R^2)$.  As $r_2\rightarrow \infty$
and $r_1\rightarrow 0$, we get
\[
  (m-1+\frac{2}{p})\norm{\frac{g}{r}}^p_{L^p(\R^2)}
  \leq \norm{\frac{g}{r}}^{p-1}_{L^p(\R^2)}
  \norm{g_r-\frac{m}{r}g}_{L^p(\R^2)},
\] 
where we used $p>2$ and the boundedness of $g$.
This completes the proof.
\end{proof}

\begin{remark}
It is essential to assume $g$ is bounded, as can be seen by the
example $g(r)=r^m$.  If we assume in the above Lemma that $g(r)=o(1)$
as $r\rightarrow\infty$, then we also have
\[
  \norm{\frac{g(r)}{r}}_{L^2(\R^2)}\leq 
  \frac 1m\norm{(\partial_r-\frac{m}{r})g(r)}_{L^2(\R^2)}.
\]
\end{remark}


Using Lemma~\ref{Lemma4-1}, we prove an $L^p$-version
of Lemma~\ref{Lemma3-2}.
\begin{lemma} 
\label{Lemma4-3} Let $2 < p <\infty$. There exists $\epsilon > 0$
such that if $f(r)$ is a radial function satisfying 
$|\lan f ,h_1 \ran _X| \leq \epsilon \|f\|_X$, then
\[
  \norm{f_r}_{L^p(\R^2)}+\norm{\frac{f}{r}}_{L^p(\R^2)}
  \leq C\big(\norm{L_0f}_{L^p(\R^2)}+\norm{L_0f}_{L^2(\R^2)}\big).
\]
Recall $L_0 f := f_r + \frac{m}{r} h_3 f $.
\end{lemma}  
\begin{proof}
We note first that it suffices to prove 
$\|f/r\|_{L^p}\leq C(\|L_0f\|_{L^p} +\|L_0f\|_{L^2})$,
since
\[
  \norm{f_r}_{L^p}\leq C(\norm{f_r+\frac{m}{r}h_3f}_{L^p}
  +\norm{\frac{m}{r}h_3f}_{L^p})
  \leq C(\norm{L_0f}_{L^p}+\norm{\frac{f}{r}}_{L^p}).
\]
Let $\varphi: [0,\infty) \rightarrow \R$ 
be a standard cut-off function with
\[
  0\leq\varphi\leq 1, \,\,\,
  \varphi(r) \equiv 1 \mbox{ for } r<1, \,\,\, \mbox{ and }\,\,\,
  \varphi(r) \equiv 0 \mbox{ for } r>2.
\]
\[
  \int_0^{\infty}|\frac{f}{r}|^p rdr
  \leq C(\int_0^{\infty}|\frac{f\varphi}{r}|^p  rdr
  +\int_0^{\infty}|\frac{f (1-\varphi)}{r}|^p rdr)
  \equiv C(I+II).
\]
We consider the second term $II$. 
Since $1-\varphi \equiv 0$ if $r<1$, we have
\[
\begin{split}
  II &= \int_0^{\infty} \abs{\frac{f(1-\varphi)}{r}}^prdr
  \leq \norm{f}^{p-2}_{L^{\infty}(\R^2)}
  \int_0^{\infty}\abs{\frac{f}{r}}^2 (1-\varphi)rdr \\
  &\leq C\norm{f}^{p-2}_{L^{\infty}(\R^2)}\norm{\frac{f}{r}}^2_{L^2(\R^2)}
  \leq C\norm{f}^{p-2}_{L^{\infty}(\R^2)}\norm{L_0f}^2_{L^2(\R^2)},
\end{split}
\]
where we used Lemma \ref{Lemma3-2}.
Next we consider the term $I$. Using Lemma~\ref{Lemma4-1}, 
we have
\[
  \norm{\frac{f\varphi}{r}}^p_{L^p}\leq C\norm{(f\varphi)_r
  -\frac{m}{r}f\varphi}^p_{L^p}
  \leq C\big(\norm{\varphi_r f}^p_{L^p}+\norm{(f_r+\frac{mh_3}{r}f)
  \varphi}^p_{L^p}
  +\norm{\frac{m}{r}(1+h_3)f\varphi}^p_{L^p}\big).
\]
Since $\varphi_r$ is supported only on $(1,2)$, and 
$(1+h_3)/r$ is bounded, we have
\[
  \norm{\frac{f\varphi}{r}}^p_{L^p}\leq 
  C(\norm{f}^p_{L^{\infty}}+\norm{L_0f}^p_{L^p}).
\]
Since, by Lemma~\ref{Lemma3-2},
$\|f\|_{L^{\infty}}\leq C\|L_0f\|_{L^2}$, we obtain
\[
  \| f/r \|^p_{L^p} \leq 
  C(\|L_0f\|^p_{L^p}+\|L_0f\|^p_{L^2}),
\]
completing the proof.
\end{proof}

\subsection{Some harmonic map estimates}

Here we prove some facts about the family
$\sO_m$ of $m$-equivariant harmonic maps.

The first lemma shows that if m-equivariant harmonic maps $h=h^{0,1}$ and
$h^{\alpha,s}$ are close in the sense of energy, then
$\alpha$ and $s$ are also close to $0$ and $1$, respectively.
   
\begin{lemma} 
\label{Lemma2-3}
Let $0<s<\infty$ and  $-\pi \leq \al < \pi$.
There exists $\epsilon > 0$ and $C > 0$ such that 
if $\int_{\R^2}\abs{\nabla \left(e^{m \th R} 
(h-h^{s,\al})\right)}^2 \,dx < \de^2$
for any $\de < \epsilon$, then
$|\al| + |s-1| \le C\de$.
\end{lemma}
\begin{proof}
Consider the case, $s>1$ (the case $s<1$ can be treated in the same way).
We first note that 
\[
  h-h^{s,\al}= \left( \begin{array}{c} 
  h_1(r)-h_1(r/s)\cos\al\\ 
  -h_1(r/s)\sin\al \\ 
  h_3(r)-h_3(r/s) \end{array} \right).
\]
Our assumption is
\[
  \int_{\R^2} \abs{\nabla \big(e^{m \th R} (h-h^{s,\al})\big)}^2 \,dx
  = 2\pi \int_0^{\infty} \big(\abs{\pd_r(h-h^{s,\al})}^2
  +\frac{m^2}{r^2}\abs{R(h-h^{s,\al})}^2\big)rdr
  < \de^2.
\]
For $0 \leq r \leq 1$, we have $(s^{2m}-1)(1-r^{2m}) \geq 0$
which, rearranged, yields 
\[
  h_1(r/s) \leq \frac{2s^m}{s^{2m}+1}h_1(r)
  \leq h_1(r), \quad 0\leq r\leq 1.
\]
Using this inequality, we find
\[
\begin{split}
  \de^2 &\geq \int_{\R^2} \abs{\nabla \big(e^{m \th R} 
  (h-h^{s,\al})\big)}^2 \,dx
  \geq 2\pi \int_0^{1}|\pd_r(h_3(r)-h_3(\frac{r}{s}))|^2rdr \\
  &= 2\pi \int_0^{1}\frac{m^2}{r^2}
  \big(h_1^2(r)-\frac{1}{s^2}h_1^2(\frac{r}{s})\big)^2rdr
  \geq 2\pi \int_0^{1}\frac{m^2}{r^2}
  \big(h_1^2(r)-h_1^2(\frac{r}{s})\big)^2rdr \\
  &\geq 2\pi m^2 \frac{(s^{2m}-1)^4}{(s^{2m}+1)^4}
  \int_0^{1} \frac{h_1^4(r)}{r^2} rdr
  = C \frac{(s^{2m}-1)^4}{(s^{2m}+1)^4}.
\end{split}
\]
It follows that $\abs{s-1}\leq C\delta^{1/2}$ 
if $\delta$ is sufficiently small.
Now $g(s) := \| \pd_r(h_3(r) - h_3(r/s)) \|_{L^2}^2$
is a smooth function of $s$ with
$g(1) = g'(1) = 0$ and $g''(1) > 0$,
and so by Taylor's theorem, we have
$g(s) \geq C(s-1)^2$ for some $C>0$, and for
$\abs{s-1} \leq C\delta^{1/2}$, $\delta$ sufficiently small.
Thus $|s-1| \leq C \delta$, as required. 

Next consider the second component of $h-h^{s,\al}$:
\[
  \de^2 > 2\pi \sin^2(\al) \int_0^{\infty}
  \abs{\pd_r h_1(\frac{r}{s})}^2rdr \geq C \sin^2(\alpha)
\]
and so $|\sin(\al)| < C\de$ for sufficiently small $\de$.
Finally, use
\[
  (1-\cos(\alpha))h_1(r) = (h - h^{s,\alpha})_1
  + \cos(\alpha)(h_1(r/s)-h_1(r))
\]
together with the previous results to arrive at
$1-\cos(\alpha) \leq C \delta$, from which (for $\alpha \in [-\pi,\pi)$)
$|\alpha| \leq C \delta$ follows.
\end{proof}

\medskip

The next lemma is a bound on the curvature
of the family $\sO_m$ of $m$-equivariant harmonic maps.
\begin{lemma}
\label{lem:curvature}
There are $\e > 0$ and $C > 0$ such that if 
\[
  \| e^{m \theta R} (h^{s_1,\al_1}(r) - h^{s_2,\al_2}(r)) \|_{\dot H^1}
  < \epsilon
\]
then setting $\bar{s} := \frac12[s_1+s_2]$,
$\bar{\al} := \frac12[\al_1+\al_2]$, and 
$\bar{h} := \frac12[h^{s_1,\al_1} + h^{s_2,\al_2}]$,
we have
\begin{equation}
\label{eq:curvebound}
  \| e^{m \theta R} (\bar{h}(r) - h^{\bar{s},\bar{\al}}(r)) \|_{\dot H^1} 
  < C \| e^{m \theta R} (h^{s_1,\al_1}(r) - h^{s_2,\al_2}(r)) \|_{\dot H^1}^2. 
\end{equation}
\end{lemma}
\begin{proof}
By rotating and rescaling, we may assume 
$(s_1,\alpha_1) = (1,0)$.
If $\epsilon$ is sufficiently small,
Lemma~\ref{Lemma2-3} gives (taking $\alpha_2 \in [-\pi,\pi)$)
\begin{equation}
\label{eq:parambound}
  |s_2-1| + |\alpha_2| \leq C
  \| e^{m \theta R} (h^{1,0}(r) - h^{s_2,\al_2}(r)) \|_{\dot H^1}
  \leq C \epsilon.
\end{equation}
Now set $s(t) := \bar{s} + (t/2)(s_2-1)$,
$\al(t) := \bar{\al} + (t/2)\al_2$,
and $\phi(t) := h^{s(t),\al(t)}$.
Then
\[
\begin{split}
  \bar{h} - h^{\bar{s},\bar{\al}} &=
  \frac{1}{2}[\phi(-1)-\phi(0) + \phi(1)-\phi(0)]
  = \frac{1}{2}[\int_0^1 \phi'(t) dt - \int_{-1}^0 \phi'(t) dt] \\
  &= \frac{1}{2} \int_0^1 [\phi'(t)-\phi'(-t)] dt
  = \frac{1}{2} \int_0^1 \int_{-t}^t \phi''(\tau) d\tau.
\end{split}
\]
Using~(\ref{eq:parambound}), we have
\[
  \| e^{m \theta R} \phi''(\tau) \|_{\dot H^1}
  \leq C[(s_2-1)^2 + \al_2^2]
  \leq C \| e^{m \theta R} (h^{1,0}(r) - h^{s_2,\al_2}(r)) \|^2_{\dot H^1}
\]
and~(\ref{eq:curvebound}) follows. 
\end{proof}

\medskip

Our next lemma gives $L^\infty$ smallness
for $\dot H^1$-small perturbations of harmonic maps.
\begin{lemma}
\label{lem:pointwise}
For $u \in \Sigma_m$, set $s = s(u)$, $\al = \al(u)$, and write
\[
  u(r,\theta) = e^{m \theta R} v(r), \quad\quad
  v(r) = e^{\alpha R}[h(r/s) + \xi(r/s)].
\]
There exists $\epsilon > 0$ and $C > 0$ such that if 
$\delta_0 < \epsilon$ and  
\begin{equation}
\label{eq:xibound2}
  \norm{\xi_r}^2_{L^2(\R^2)} + \norm{\frac{m}{r} R\xi}^2_{L^2(\R^2)} =
  \| u-e^{m\theta R}h^{s,\al} \|_{\dot H^1}^2
  < \delta^2_0,
\end{equation}
then
\[
  \| \xi \|_{L^\infty} \leq C \delta_0.
\]
\end{lemma}
\begin{proof}
Without loss of generality, we may assume $s=1$ and $\alpha=0$.
It follows immediately from~(\ref{eq:xibound2}) and 
Lemma~\ref{Lemma2-1} that
\[
  \|\xi_i\|_{\infty} \leq C\delta_0, \quad\quad i=1,2.
\]
For $\xi_3$, we have, as yet, only 
$\norm{(\xi_3)_r}_{L^2} < \delta_0$,
and so our aim is to show that 
$\norm{\xi_3}_{L^{\infty}}\leq C\delta_0$.
Under our change of variable,
$\tilde{\xi}(y) := \xi(m\log(r))$,
it suffices to prove that 
$\norm{\tilde{\xi}_3}_{L^2(\R)}\leq C\delta_0$, since 
$\norm{\tilde{\xi}_3}_{L^{\infty}(\R)}\leq C \norm{\tilde{\xi}_3}_{H^1(\R)}$.
By the continuity and boundary conditions of $v(r)$ (for $u \in \Sigma_m$),
there must exist $y_0 \in \R$ such that $\tilde{v}_3(y_0)=0$.
Note that since $\tilde{h}_3(y) = \tanh(y)$, we have
\begin{equation}
\label{eq:estimate}
  |\tilde{v}_3^2(y) - \tanh^2(y)| 
  = |\sum_{j=1}^2 (\tilde{v}_j^2 - \tilde{h}_j^2)|
  \leq C (|\xi_1| + |\xi_2|)
\end{equation}
and in particular 
\[
  \tanh^2(y_0) \leq 
  C(\| \xi_1 \|_{L^\infty} + \| \xi_2 \|_{L^\infty})
  < C \delta_0.
\]
So for $-1 \leq y \leq 1$,
\begin{equation}
\label{eq:localbit}
\begin{split}
  |\tilde{v}_3(y) - \tanh(y)|
  &= |\int_{y_0}^y (\tilde{v}_3' - \tanh'(y)) dy - \tanh(y_0)| \\
  &\leq C (\|\tilde{v}_3' - \tanh'\|_{L^2} + \delta_0^{1/2})
  \leq C \delta_0^{1/2},
\end{split}  
\end{equation}
and in particular, for $\delta_0$ sufficiently small,
$|\tilde{v}_3(\pm 1)| > (1/2)\tanh(1)$.
Then with the aid of Lemma~\ref{Lemma3-4}, for
$\delta_0$ sufficiently small, we have 
$|\tilde{v}_3(y)| > (1/4)\tanh(1)$ for $|y| \geq 1$.
Estimate~(\ref{eq:estimate}) then yields
\[
  \int_{|y| \geq 1} (\tilde{v}_3(y) - \tanh(y))^2 dy
  \leq C \delta_0^2.
\]
which also gives us 
$\sup_{|y| \geq 1} |\tilde{v}_3(y) - \tanh(y)| \leq C \delta_0$
(Sobolev embedding), and in particular
$|\tilde{v}_3(1) - \tanh(1)| \leq C \delta_0$.
Finally, we get the same result for $|y| < 1$ by integrating
the derivative:
\[
  |\tilde{v}_3(y) - \tanh(y)| \leq
  |\int_y^1 ( \tilde{v}_3' - \tanh' ) dy| + C \delta_0
  < C \delta_0.
\]
This completes the proof.
\end{proof}

\subsection{Perturbation is bounded by $q$}

Here we prove estimates used in Section~\ref{sec:flow}.
We show that $z(r)$
is controlled by $q(r)$, where, recall,
\[
  v_r - \frac{m}{r}J^v Rv
  = q_1 \hat{e} + q_2 J^v \hat{e};
  \quad\quad D_r^v \hat{e} \equiv 0
\]
and
\begin{equation}
\label{eq:v}
  v(r) = e^{\al R}[h(r/s) + \xi(r/s)], \quad\quad
  \xi = z_1 e + z_2 J^h e + \gamma h.
\end{equation}
\begin{lemma}
\label{Lemma4-5}
Let $2\leq p<\infty$. For $\delta$ sufficiently small,
\begin{equation}
  \norm{z_r}_{L^p(\R^2)} + \norm{\frac{z}{r}}_{L^p(\R^2)} 
  \leq C \bke{ s^{1-2/p} \norm{q}_{L^p(\R^2)}
  + \norm{q}_{L^2(\R^2)} },
\label{lemma4-5:10}
\end{equation}
\begin{equation}
  \norm{z_{rr}}_{L^2(\R^2)} 
  \le C \bke{ s\norm{q_r}_{L^2(\R^2)} + s\norm{\frac{q}{r}}_{L^2(\R^2)}
  + s\norm{q}^2_{L^4(\R^2)} + \| q \|_{L^2(\R^2)} },
\label{Lemma4-6:10}
\end{equation}
and
\begin{equation}
  \| u \|_{\dot H^2(\R^2)} \le C \bke{ 
  \frac{1}{s} + \norm{q_r}_{L^2(\R^2)} + \norm{\frac{q}{r}}_{L^2(\R^2)}
  + \norm{q}^2_{L^4(\R^2)} + \frac{1}{s}\| q \|_{L^2(\R^2)} }.
\label{eq:ubound}
\end{equation}

\end{lemma}
\begin{proof}
We will first show the following:
\begin{equation}
  \norm{z_r}_{L^2(\R^2)} + \norm{\frac{z}{r}}_{L^2(\R^2)} 
  \le C \norm{q}_{L^2(\R^2)}.
\label{Lemma4-5:20}
\end{equation}
By~(\ref{eq:v}), we have
\begin{equation}
  e^{-\al R} s(q\hat{e})(sr) =
  (L_0z)e + (\gamma h)_r+\frac{2m}{r}h_3\gamma h
  +\frac{m}{r}\xi_3\xi.
\label{Lemma4-5:30}
\end{equation}
Since $\|z\|_X\leq C\|L_0z\|_{L^2}$ by~(\ref{eq:star}),
it suffices to prove that $\|L_0z\|\leq C\|q\|_{L^2}$.
We first show that $\|\xi_r\|_{L^2}+\|\xi/r\|_{L^2}\leq C\|z\|_X$.
Indeed, since
$(J^h e)_r = -(m/r)h_1h$ and $h_r=(m/r)h_1J^he$, we find
\[
  \xi_r =  z_r e
  -z_2\frac{m}{r} h_1h +\gamma_r h
  +\frac{m}{r} \gamma h_1 J^h e.
\]
Therefore, since $\gamma = O(|z|^2)$ and 
$\gamma_r = O(|z||z_r|)$, 
we obtain
\[
  \norm{\xi_r}_{L^2}\leq \norm{z_r}_{L^2}+\norm{\frac{z}{r}}_{L^2}+
  \norm{z}_{L^{\infty}}\norm{z_r}_{L^2}+\norm{z}_{L^{\infty}}\norm{\frac{z}{r}}_{L^2},
\]
where we used the boundedness of $h$.
By~(\ref{eq:smallness}), we have 
$\|z\|_{L^{\infty}} \leq C\de$, 
which can be chosen sufficiently small to yield
$\|\xi_r\|_{L^2} \leq C\|z\|_X$.
In a similar manner, we can show 
$\|\xi/r\|_{L^2} \leq C \|z/r\|_{L^2} \leq C \|z\|_X$. 
Combining, we obtain 
\begin{equation}
\label{eq:one}
  \| \xi_r\|_{L^2} + \| \xi/r \|_{L^2} \leq C \|z\|_X.
\end{equation}
Now we are ready to prove $\|L_0z\|_{L^2}\leq C\|q\|_{L^2}$.
Using again $\gamma = O(|z|^2)$,  
$\gamma_r = O(|z||z_r|)$, and the boundedness of $h$, we find 
\[
\begin{split}
  \norm{(\gamma h)_r+\frac{2m}{r}h_3\gamma h+\frac{m}{r}\xi_3\xi}_{L^2}
  &\leq C\norm{z}_{L^{\infty}}\big(\norm{z_r}_{L^2}+\norm{\frac{z}{r}}_{L^2}\big)
  +C\norm{\xi}_{L^{\infty}}\norm{\frac{\xi}{r}}_{L^2} \\
  &\leq C\norm{\xi}_{L^{\infty}}\big(\norm{z}_X+\norm{\frac{\xi}{r}}_{L^2}\big)
  \leq C\norm{\xi}_{L^{\infty}}\norm{z}_X \\
  &\leq C\norm{\xi}_{L^{\infty}}\norm{L_0z}_{L^2},
\end{split}
\]
where we used~(\ref{eq:one}) and
$\|z\|_{L^{\infty}} \leq \|\xi\|_{L^{\infty}}$.
Thus we have
\[
  \norm{L_0z}_{L^2} \leq \norm{s q(s \cdot)}_{L^2} +
  \norm{(\gamma h)_r+\frac{2m}{r}h_3\gamma h+\frac{m}{r}\xi_3\xi}_{L^2}
  \leq \norm{q}_{L^2}+C\norm{\xi}_{L^{\infty}}
  \norm{L_0z}_{L^2}.
\]
Since 
$\|\xi\|_{L^{\infty}} \leq C\|\xi\|_X \leq C\|z\|_X \leq C\delta$ 
can be taken sufficiently small, the above inequality implies
\begin{equation}
\label{eq:two}
  \|L_0z\|_{L^2}\leq C\|q\|_{L^2},
\end{equation}
which completes the proof of \eqref{Lemma4-5:20}.
Similarly, for any $p$ with $2\leq p<\infty$, we have
\[
\begin{split}
  \norm{L_0z}_{L^p} &\leq 
  C\bke{ \norm{sq(s \cdot)}_{L^p}+\norm{z}_{L^{\infty}}
  \big(\norm{z_r}_{L^p}+\norm{\frac{z}{r}}_{L^p}\big)} \\
  & \leq C\bke{s^{1-2/p}\norm{q}_{L^p}+\norm{z}_{L^{\infty}}
   \big(\norm{L_0z}_{L^p} + \norm{L_0 z}_{L^2}\big)},
\end{split}
\]
where we used Lemma~\ref{Lemma4-3}.
Since $\|z\|_{L^{\infty}}$ can be taken sufficiently small, 
and using~(\ref{eq:two}), we finally have 
$\|L_0z\|_{L^p}\leq C(s^{1-2/p}\|q\|_{L^p} + \|q\|_{L^2})$,
completing the proof of~(\ref{lemma4-5:10}).

Next we prove~(\ref{Lemma4-6:10}).
We first show
\begin{equation}
  \norm{z_{rr}}_{L^2}
  \leq C \left(\norm{\pd_r L(ze)}_{L^2}+\norm{\frac{L(ze)}{r}}_{L^2}
  + \norm{L_0z}_{L^2}\right).
\label{Lemma4-6:20}
\end{equation}
Indeed, recalling $L(ze)=z_re+(m/r)h_3ze$, we have
\[
\begin{split}
  (\pd_r+\frac{1}{r}+\frac{m}{r}h_3)z_re &= 
(\pd_r+\frac{1}{r})z_re+\frac{m}{r}h_3z_re\\
  &=(\pd_r+\frac{1}{r})(L(ze)-\frac{m}{r}h_3ze)
+\frac{m}{r}h_3z_re\\
&=(\pd_r+\frac{1}{r})L(ze)-\frac{m^2}{r^2}h^2_1ze
+\frac{m^2}{r^2}z_2h_3h_1 h.
\end{split}
\]
Set 
\[
  \eta := e^{m\theta R}ze 
  \quad \mbox{ and } \quad
  H := e^{m\theta R}h.
\]
Since $\pd_j \eta=D^{H}_j \eta -(\pd_j H\cdot \eta)H$ 
and $(h_3)_r = (m/r)h^2_1$, we have 
\[
\begin{split}
  e^{-m\th R}D^{H}_jD^{H}_j\eta 
&=\bke{(\pd_r+\frac{1}{r})z_r}e-\frac{m^2}{r^2}h^2_3ze\\
&=(\pd_r+\frac{1}{r})(z_re)+\frac{m}{r}z'_2h_1h-\frac{m^2}{r^2}h^2_3ze\\
&=(\pd_r+\frac{1}{r}+\frac{m}{r}h_3)(z_re)-
\frac{m}{r}h_3L(ze)+\frac{m}{r}z'_2h_1h,
\end{split}
\]
and therefore we obtain 
\[
\begin{split}
  \norm{D^{H}_jD^{H}_j\eta}_{L^2}
  &=\norm{\bke{(\pd_r+\frac{1}{r})z_r}e-\frac{m^2}{r^2}h^2_3ze}_{L^2}\\
&=\norm{(\pd_r+\frac{1}{r})L(ze)-\frac{m}{r}h_3L(ze)
-\frac{m^2}{r^2}h^2_1ze
+\frac{m^2}{r^2}z_2h_3h_1 h
+\frac{m}{r}z'_2h_1h}_{L^2}\\
&\leq C\bke{\norm{(\pd_r+\frac{1}{r})L(ze)}_{L^2}+\norm{\frac{L(ze)}{r}}_{L^2}
  + \norm{z'}_{L^2}+\norm{\frac{z}{r}}_{L^2}} \\
  &\leq C\bke{\norm{\pd_rL(ze)}_{L^2}+\norm{\frac{L(ze)}{r}}_{L^2}
  + \norm{L_0z}_{L^2}},
\end{split}
\]
where we used $|h_1/r|$ bounded.
Since $\pd_j \eta=D^{H}_j \eta -(\pd_j H\cdot \eta)H$, we have
\[
\begin{split}
  \Delta \eta &= \pd_j\pd_j \eta
  =\pd_j\bke{D^{H}_j \eta-(\pd_j H\cdot \eta)H}
  =\pd_jD^{H}_j \eta-\pd_j\bke{(\pd_j H\cdot \eta)H} \\
  &=D^{H}_jD^{H}_j \eta-(\pd_j H\cdot D^{H}_j \eta)H
   -(\pd_j\pd_j H\cdot \eta)H-(\pd_j H\cdot \pd_j \eta)H
   -(\pd_j H\cdot \eta)\pd_j H.
\end{split}
\]
Therefore, we obtain
\[
\begin{split}
  \norm{\Delta \eta}_{L^2} &\leq \norm{D^{H}_jD^{H}_j \eta}_{L^2}
   +\norm{(\pd_j H\cdot D^{H}_j \eta)H}_{L^2}
   +\norm{(\pd_j H\cdot \pd_j \eta)H}_{L^2}\\
  &\quad +\norm{(\pd_j\pd_j H\cdot \eta)H}_{L^2}
         +\norm{(\pd_j H\cdot \eta)\pd_j H}_{L^2}\\
  &\leq \norm{D^{H}_jD^{H}_j \eta}_{L^2}
  +C\bke{\norm{\nabla \eta}_{L^2}
  +\norm{\frac{\eta}{r}}_{L^2}} \\
  &\leq \norm{D^{H}_jD^{H}_j \eta}_{L^2}+C\bke{\norm{z_r}_{L^2}
  +\norm{\frac{z}{r}}_{L^2}}.
\end{split}
\]
Thus
\begin{equation}
\label{eq:eta}
  \norm{\eta}_{\dot{H}^2}\leq C\norm{\Delta \eta}_{L^2}
  \leq C\bke{\norm{\pd_rL(ze)}_{L^2}+\norm{\frac{L(ze)}{r}}_{L^2}
  + \norm{L_0z}_{L^2}}.
\end{equation}
Let $A=(ze)_{rr}-(ze)_r/r+m^2ze/r^2, B=(ze)_{rr},$ and 
$E=(mR(ze)_r)/r-mRze/r^2$. Then direct calculations show 
\[
(e^{m\theta R}ze)_{xx}=e^{m\theta R}
\bke{-\frac{y^2}{r^2}A+B-2\frac{xy}{r^2}E},
\]
\[
(e^{m\theta R}ze)_{yy}=e^{m\theta R}
\bke{-\frac{x^2}{r^2}A+B+2\frac{xy}{r^2}E},
\]
\[
(e^{m\theta R}ze)_{xy}=e^{m\theta R}
\bke{\frac{xy}{r^2}A
+(\frac{x^2}{r^2}-\frac{y^2}{r^2})E}.
\]
Each of $A, B,$ and $E$ can be expressed in terms of 
combinations of second derivatives of $\eta = e^{m\theta R}ze$.
This implies, in particular, the estimate~\eqref{Lemma4-6:20}.

It remains to control
$\norm{\pd_r L(ze)}_{L^2}+\norm{L(ze)/r}_{L^2}$.
Noting that $\hat{e}_r=-(v_r\cdot\hat{e})v$,
$(J^v\hat{e})_r=(v_r \cdot \hat{e})v$, so
\[
  (q\hat{e})_r = q_r\hat{e} - q_1(v_r\cdot\hat{e})v
  + q_2(v_r\cdot\hat{e})v,
\]
and recalling $v_r(r) = e^{\al R}\frac{1}{s}(h_r(r/s)+\xi_r(r/s))$,
we have
\[
\begin{split}
  \norm{\pd_r(q\hat{e})}_{L^2} & \leq 
C\bke{\norm{q_r}_{L^2}+\norm{\frac{q}{r}}_{L^2}
  +\norm{q \frac{1}{s}z_r(\cdot/s)}_{L^2}} \\
  & \leq C\bke{\norm{q_r}_{L^2}+\norm{\frac{q}{r}}_{L^2}
  +\norm{q}^2_{L^4} + s^{-1}\norm{q}_{L^2}^2 },
\end{split}
\]
where we used~\eqref{lemma4-5:10}.
Taking the derivative of Equation~\eqref{Lemma4-5:30}, we get
\[
\begin{split}
  \norm{\pd_r L(ze)}_{L^2} & \leq C \bke{s\norm{q_r}_{L^2}
  + s\norm{\frac{q}{r}}_{L^2} + s\norm{q}_{L^4}^2 + \norm{q}_{L^2}^2 \right. \\
  & \quad \quad \left. +\norm{(\gamma h)_{rr}+(2\frac{m}{r}h_3\gamma h)_r
  +(\frac{m}{r}\xi_3\xi)_r}_{L^2} }.
\end{split}
\]
We consider first $(\gamma h)_{rr}$. Using 
$|h_r| + |rh_{rr}| <C$, we have
\[
\begin{split}
  \norm{(\gamma h)_{rr}}_{L^2}
  &\leq C \norm{|\gamma_{rr}|+|\gamma_rh_r|+|\gamma h_{rr}|}_{L^2} \\
  &\leq C\bke{ \norm{z}_{L^{\infty}}\norm{z_{rr}}_{L^2} + \|z_r\|^2_{L^4}
  + \norm{|z_r|+|\frac{z}{r}|}_{L^2} } \\
  &\leq C \bke{ \norm{z}_{L^{\infty}}\norm{z_{rr}}_{L^2}
  + s\|q\|^2_{L^4} + \norm{q}^2_{L^2} + \norm{q}_{L^2} }.
\end{split}
\]
Next we consider $(\frac{m}{r}h_3\gamma h)_r$. 
In a similar manner, we find
\[
\begin{split}
  \norm{(\frac{m}{r}h_3\gamma h)_r}_{L^2} &\leq 
  C \norm{ |\frac{\gamma}{r^2}|+ |\frac{\gamma_r}{r}|}_{L^2}
  \leq C\bke{ \norm{\frac{z}{r}}^2_{L^4}
           +\norm{z_r}_{L^4}\norm{\frac{z}{r}}_{L^4} } \\
  &\leq C\bke{s\|q\|^2_{L^4} + \|q\|^2_{L^2}}.
\end{split}
\]
For the term $(\xi_3\xi/r)_r$, we have the estimate:
\[
  \norm{(\frac{m}{r}\xi_3\xi)_r}_{L^2}
  \leq C \bke{s\|q\|^2_{L^4} + \norm{q}^2_{L^2}}.
\]
Following a similar procedure for $\|L(ze)/r\|_{L^2}$,
we obtain
\[
\begin{split}
  \norm{\pd_r L(ze)}_{L^2} + \norm{\frac{L(ze)}{r}}_{L^2}
  &\leq C\bke{ s\norm{q_r}_{L^2} + s\norm{\frac{q}{r}}_{L^2}
  +\norm{z}_{L^{\infty}}\norm{z_{rr}}_{L^2} + s\|q\|^2_{L^4} 
  + \norm{q}^2_{L^2} } \\
  &\leq C\big( s\norm{q_r}_{L^2} + s\norm{\frac{q}{r}}_{L^2}
  +\norm{z}_{L^{\infty}}
   \bke{ \norm{\pd_r L(ze)}_{L^2}+\norm{\frac{L(ze)}{r}}_{L^2} } \\
  & \quad\quad +s\|q\|^2_{L^4} + \norm{q}^2_{L^2} \big),
\end{split}
\]
using~\eqref{Lemma4-6:20}. 
Since $\norm{z}_{L^{\infty}}$ can be taken sufficiently small, 
we conclude 
\begin{equation}
\label{eq:lastbound}
  \norm{\pd_r L(ze)}_{L^2}+\norm{\frac{L(ze)}{r}}_{L^2}
  \leq C\bke{ s\norm{q_r}_{L^2} + s\norm{\frac{q}{r}}_{L^2}
  +s\|q\|^2_{L^4} + \norm{q}_{L^2}},
\end{equation}
having used the smallness of $\| q \|_{L^2}$.
Combining~(\ref{eq:lastbound}) 
with~(\ref{Lemma4-6:20}) completes the proof of~(\ref{Lemma4-6:10}).

It remains to prove~(\ref{eq:ubound}).
Since $u(x) = e^{(m \theta + \al) R}\bke{h(r/s)+\xi(r/s)}$, 
it is straightforward to check that for $\delta$ sufficiently small,
\[
  \| u \|_{\dot H^2} \leq \frac{C}{s} \left(
  1 + \| \eta \|_{\dot H^2} \right)
\]
and so~(\ref{eq:ubound}) follows from~(\ref{eq:eta})
and~(\ref{eq:lastbound}).
This completes the proof of Lemma~\ref{Lemma4-5}.
\end{proof}

\section*{Acknowledgements}

We thank Nai-Heng Chang, Jalal Shatah, and Chongchun Zeng for sharing
with us their results on self-similar singular solutions, which are
not used in this paper, but inspired us nonetheless. The
first and third authors are partially supported by NSERC grants nos.
22R80976 and 22R81253.  The second author is partially supported by a
PIMS PDF.


\end{document}